\documentclass[11pt,english]{article}
\usepackage[T1]{fontenc}
\usepackage{color}
\usepackage{float}
\usepackage{amsmath}
\usepackage{amssymb}
\usepackage{graphicx}
\usepackage{setspace}
\usepackage[authoryear]{natbib}
\makeatletter
\usepackage{amsmath}
\usepackage{sectsty} 
\usepackage{pslatex}

\sectionfont{\sffamily\large}
\subsectionfont{\sffamily\normalsize}
\allowdisplaybreaks

\newcounter{eqnum}[section] 
\newtheorem{proposition}{Proposition} 
 
\newtheorem{lemma}{Lemma} 
\newtheorem{corollary}{Corollary}

\newcommand{\eproof}{\mbox{}\hfill{\rule{8pt}{8pt}}}

\makeatother

\usepackage{babel}
\begin{document}

\title{\textbf{Invest or Exit? Optimal Decisions in the Face of a Declining
Profit Stream}\thanks{DOI 10.1287/opre.1090.0740}\textbf{ }}

\author{H. Dharma Kwon{\normalsize{} }\thanks{Department of Business Administration, University of Illinois at Urbana-Champaign,
Champaign, IL 61820.\protect \\
dhkwon@illinois.edu}}

\date{{\normalsize{}April 24, 2009 }}
\maketitle
\begin{abstract}
Even in the face of deteriorating and highly volatile demand, firms
often invest in, rather than discard, aging technologies. In order
to study this phenomenon, we model the firm's profit stream as a Brownian
motion with negative drift. At each point in time, the firm can continue
operations, or it can stop and exit the project. In addition, there
is a one-time option to make an investment which boosts the project's
profit rate. Using stochastic analysis, we show that the optimal policy
always exists and that it is characterized by three thresholds. There
are investment and exit thresholds before investment, and there is
a threshold for exit after investment. We also effect a comparative
statics analysis of the thresholds with respect to the drift and the
volatility of the Brownian motion. When the profit boost upon investment
is sufficiently large, we find a novel result:  the investment threshold
decreases in volatility.
\end{abstract}

\section{Introduction\label{sec:introduction}}

The computer disk drive industry underwent a series of disruptive
architectural innovations (\citealt{Christensen1992}). Until the
mid-1970's, 14-inch hard disk drives dominated the mainframe computer
disk drive market. Between 1978 and 1980, several new entrants introduced
8-inch disk drives which were initially sold to minicomputer manufacturers
because their recording capacity was too small and the cost per megabyte
was too high for mainframe computers. As the performance of 8-inch
drives kept improving, the entrants quickly encroached upon the mainframe
computer market. By the mid-1980's, 8-inch drives dominated the mainframe
market and rendered 14-inch drives obsolete. Nevertheless, among the
dozen or so established manufacturers of 14-inch drives, two thirds
of them never introduced 8-inch drives. Instead, they continued to
enhance the recording capacity of the extant 14-inch drives in order
to appeal to the higher end mainframe market (\citealt{Christensen2000},
p. 19). Eventually, all 14-inch drive manufacturers, except those
that were vertically integrated, were forced out of the disk drive
market. This pattern of industry-wide disruption emanating from the
introduction of a successful new technology is a commonplace rather
than an isolated incident; as such, it deserves serious attention.
Even  8-inch drives were eventually superseded by 5.25-inch drives.
Currently, the computer disk drive industry is in the process of yet
another architectural transition, one from hard disk drives to flash
solid state disks (used in USB stick drives). 

This paper focuses upon the difficult investment and exit decisions
of a firm facing a declining profit stream. With the onslaught of
disruptive technological innovation, as in the example of the disk
drive industry, a firm employing an extant technology faces a deteriorating
profit stream due to declining demand and/or prices. Faced with a
profit stream that has eroded, it might be optimal for the firm to
cease operations and avoid recurring losses. On the other hand, if
the erosion has not been too large, then it can be optimal for the
firm to make an additional investment in the project. The pressing
question is when, if ever, to invest and when to exit. Exit ought
to occur when the current profit rate is sufficiently negative; a
negative value of the profit rate, however, is not a sufficient condition
to induce exit because the option to cease operations sometime in
the future must be taken into account. Likewise, a firm must invest
in its operations in a timely fashion before the desirable investment
opportunity vanishes. In a highly volatile environment such as in
the disk drive industry, however, it is difficult to calculate the
optimal time to invest or exit because of the uncertainty in the future
demand. After we obtain the optimal policy, we will examine how increases
in uncertainty affect the optimal policy. 

In light of the declining demand, it seems counter-intuitive to invest
in the current operation. However, in the example of the computer
hard disk drive industry, the manufacturers of 14-inch disk drives
continued investment even though they faced a deteriorating profit
stream and, as it turned out, eventual displacement from the industry.
\citet{Christensen2000} finds such examples in the mechanical excavator
industry and the steel mill industry as well. \citet{Rosenberg1976}
also notes the major investments in obsolescent technologies during
the transition from the wooden sailing ship to the iron-hull steamship
and from the steam  to the diesel engine, to name just two amongst
many such examples.

The two$\:$\textcolor{black}{{} salient }features of our model are
the possibility of exit and the declining stochastic profit stream.
In particular, the firm can exit at any point in time, and we model
the firm's uncertain profit stream as a Brownian motion $X_{t}$ with
drift $\mu$ and volatility $\sigma$ where both $\mu$ and $\sigma$
are time-independent constants known to the firm. Of course, the drift
$\mu$ is the average rate of change in the profit rate, and the volatility
$\sigma$ measures the underlying uncertainty. In particular, we give
special attention to the case in which $\mu$ is negative. With this
representation, the firm's cumulative profit is the time-integral
of the Brownian motion. The firm's investment and exit decisions are
stopping times for the Brownian motion; we utilize the well-known
machinery of stochastic analysis (\citealt{Alvarez2001a}, \citealt{Borodin2002},
\citealt{Harrison1985}, \citealt{Oksendal2003}, \citealt{Peskir2006})
 to prove the existence of the optimal policy and find the optimal
stopping times. 

In Sec. \ref{sec:basic-model}, we present the basic model in which
investment is not possible. At each point in time, the firm must decide
whether to continue operations or irrevocably exit the project. The
firm seeks to maximize its expected discounted cumulative profit by
selecting the optimal time $\tau$ at which to exit, where $\tau$
is a stopping time for the Brownian motion. Utilizing results obtained
by \citet{Alvarez2001a}, we show that the optimal policy is a threshold
rule: it is optimal to continue operations until the profit rate $X_{t}$
falls below a critical threshold $\xi_{0}$, at which time it is optimal
to exit. The closed-form expression for $\xi_{0}$ is a decreasing
function of $\mu$ and $\sigma$, and it reveals that $\xi_{0}$ is
negative. 

In Sec. \ref{sec:Optimal-investment}, we extend the basic model to
include a one-time opportunity to invest in improving the extant technology:
at each point in time, the firm can (1) continue operations, (2) stop
and irrevocably exit the project, or (3) invest in the operations.
The investment increases both the current profit rate and the drift
of the profit stream by known quantities. In view of the investment
opportunity, the firm's policy is specified by three stopping times.
The firm must specify when to exit and when to invest while the investment
option is still available. If the firm already has made the investment,
then the firm must decide when to exit. Each stopping time is characterized
by a threshold. If investment has not been made, it is optimal to
exit whenever the profit rate falls below a threshold $\xi_{E}$,
and it is optimal to invest if the profit rate rises above a second
threshold $\xi_{I}$. When the current profit rate  is between $\xi_{E}$
and $\xi_{I}$, it is optimal to maintain the \emph{status quo}: continue
operations but do not invest. Because there is only one opportunity
to invest, after investment, the firm's decision problem reduces to
that of the basic model, albeit with different drift: after investment,
the firm exits when the current profit rate drops below a third threshold
$\xi_{1}$. 

After finding the optimal policy, we effect a comparative statics
analysis of the thresholds $\xi_{I}$ and $\xi_{E}$ with respect
to $\mu$ and $\sigma$. Although it is intuitively clear that the
optimal policy is characterized by thresholds, the comparative statics
analysis is neither intuitively obvious nor straightforward. In order
to obtain $\xi_{I}$ and $\xi_{E}$, we first need to solve an optimal
stopping time problem with a reward which depends on the return from
investment. The complication is that the return from investment in
turn depends on both $\mu$ and $\sigma$ because the firm will continue
operations prior to eventual exit. Nevertheless, we have been able
to effect a comparative statics analysis using a power-series expansion
method without resorting to a numerical analysis. 

We must proceed cautiously in applying  intuition from real options
theory to the comparative statics of the threshold for investment
($\xi_{I}$). For example, consider real options models in which the
return from investment is independent of $\sigma$; real options theory
has shown that, under certain mild conditions, it is optimal to wait
longer before making an irrevocable investment if the volatility of
the underlying asset increases (\citealt{Dixit1992}, \citealt{Alvarez2003}).
Waiting and observing the evolution of the value of the asset enables
the investor to avoid the downturn risk and take advantage of the
upturn potential. In accord with this intuition, we anticipate that
$\xi_{I}$ increases in $\sigma$ because the upturn potential of
the profit stream increases in $\sigma$. Indeed, if the boost in
the profit rate upon investment is small enough, then $\xi_{I}$ increases
in $\sigma$ as expected. Surprisingly, if the boost is sufficiently
large, then $\xi_{I}$ decreases in $\sigma$; this novel comparative
statics of $\xi_{I}$ results because the return from investment increases
in $\sigma$ due to the embedded option to exit, and consequently,
the investment option is more attractive with higher $\sigma$. The
ever-present exit option is the salient feature which sets our model
apart from the conventional real options models. In the operations
context, our novel comparative statics result offers cautionary advice
against blindly following the intuition inherited from real options
theory. For example, \citet{Bollen1999} shows that if the  product
life cycle (demand dynamics) is ignored, then the conventional real-option
technique tends to undervalue capacity contraction and overvalue capacity
expansion. 

This paper is organized as follows. We review related literature in
Sec. \ref{sec:related-literature}, and we present our basic model
and  some results about the basic model in Sec. \ref{sec:basic-model}.
The basic model is extended to include one investment opportunity
in Sec. \ref{sec:Optimal-investment}, and we effect the comparative
statics of the optimal policy in Sec. \ref{sec:Comparative-Statics}.
Summary of the paper is given in Sec. \ref{sec:summary}. 

\section{Related Literature\label{sec:related-literature}}

There is a rich literature on technology and process adoption. (See,
for example, \citealt{Bridges1991} for a review.) In an early paper
which formulates technology adoption as an investment problem, \citet{Barzel1968}\textsf{\textbf{\large{}
}}uses the net-present-value approach to obtain the optimal timing
of a one-time investment in the adoption of technology when the future
profit stream is deterministic. In the context of process improvement,
\citet{Porteus1985} uses the EOQ model to examine the economic trade-offs
between the cost of investment which reduces the setup cost and the
benefit from the reduced setup cost: the optimal policy is to invest
if and only if the sales rate is above a threshold. \citet{Porteus1986}
extends this work by examining a model in which lower setup costs
lead to improved quality control (lower defect rate). 

An objective of the current paper is to obtain the optimal investment
and exit policy under uncertainty. Many papers have modeled technology
adoption as a stopping time problem. (See \citealt{Hoppe2002} for
a survey of literature.) For example, \citet{Balcer1984} study the
optimal time to adopt the best currently available technology when
multiple adoptions are allowed. In their model, the \emph{timing}
and the \emph{value} of future innovations is uncertain although the
profitability of the currently available technology is known. They
show that it is optimal to adopt the best currently available technology
if the technological lag exceeds a threshold which depends upon the
multi-dimensional state: the elapsed time since last innovation and
the pace (rapidity) of technological progress. There is also a substantial
literature on Bayesian models of investment and exit; see \citet{Jensen1982},
\citet{McCardle1985}, \citet{Ryan2003}, and \citet{Ryan2005}.

One important contribution of our paper is the results regarding the
impact of uncertainty on the investment and exit decisions. \citeauthor{Dixit1992}
(1992, p. 108) points out that, as uncertainty increases, it is optimal
to wait longer before investment if (1) the investment is irreversible,
(2) the uncertainty regarding the investment is being resolved gradually
in time, and (3) the investment can be flexibly postponed. In this
vein, \citet{McDonald1986} study investment in an asset whose value
and price evolve as geometric Brownian motion. They find that the
optimal policy is a threshold rule with respect to the ratio of the
value to the price of the asset. Moreover, the investment threshold
increases in the volatility: it is optimal to postpone investment
longer as the uncertainty increases. 

A number of papers address the effect of uncertainty on technology
adoption using the real options approach. Essentially, they confirm
the conventional intuition regarding the value of waiting. \citet{Farzin1998}
study the optimal time to irreversibly switch to a new technology
when the value and the arrival date of future improvements are uncertain.
In their model, the improvement in the value of the currently available
technology follows a compound Poisson process. They allow multiple
investments in technology; again, the optimal policy is a threshold
rule. In particular, they find that the pace of adoption is slower
with the real-option method than with the suboptimal net-present-value
method. \citet{Alvarez2001} also use the real options approach to
study the optimal time to adopt a technology with an opportunity for
improvement after adoption. Once the firm adopts the technology, it
receives a revenue stream which evolves stochastically over time:
at an exponential time, an improved technology becomes available to
the firm. They show that increased market uncertainty (volatility)
increases the real-option value of adopting the initial technology. 

The real options method has also been applied to exit decisions in
a duopoly game when the profit stream is stochastic. \citet{Fine1986}
find a Nash equilibrium in stopping times in their discrete-time duopoly
game of exit from a market with declining stochastic demand. \citet{Murto2004}
studies a similar duopoly exit game in an industry in which the declining
demand follows a geometric Brownian motion; he obtains Markov-perfect
equilibria. Although these two papers analyze a duopoly model, they
also consider the exit problem of a monopolist which is similar to
our basic model. Their focus, however, is on the strategic interaction
rather than on the uncertainty.

Mathematical characteristics of general real options models have been
extensively studied in the stochastic control theory literature. In
particular, many real options models can be formulated as optimal
stopping problems. The solution method of optimal stopping time problems
consists of finding a candidate solution by solving a differential
equation and applying a verification theorem which includes smooth-pasting
conditions (that the optimal return function is continuously differentiable);
for example, see \citet{Dayanik2003}, Chapter 10 of \citet{Oksendal2003},
and Chapter IV of \citet{Peskir2006}. \citet{Alvarez2001a} considers
a class of optimal stopping problems with a profit rate and a salvage
value which are functions of a linear diffusion process. He finds
a representation of the optimal return function, a set of necessary
conditions for the optimal solution, and also conditions under which
the necessary conditions are sufficient. \citet{Alvarez2003} studies
a class of optimal stopping problems which often occur in economic
decision problems and finds general comparative statics with respect
to volatility. In particular, he characterizes a set of conditions
under which the value function and the optimal continuation region
increase in volatility. 

The solution methods for stopping time problems have been applied
to study mathematical properties of financial options (see, for example,
\citealt{Guo2001} and \citealt{Ekstrom2004}) as well as a wide variety
of real options models. For example, \citet{Wang2005} studies an
optimal stopping time problem in which the decision-maker has a sequence
of projects. For each project, there is a forced exit event at which
time the decision-maker is forced to stop the project. The forced
exits occur as a Poisson process, and the decision variable is the
sequence of entry times. If the entry cost is large enough, then the
presence of the forced exits enlarge the continuation region. The
stopping time framework can be also applied to study real options
problems under incomplete information. \citet{Decamps2005} study
the optimal timing of investment in a noisy asset whose true underlying
value (the drift of the Brownian motion) is unknown. The information
about the value of the drift is Bayesian-updated based on the observed
value of the Brownian motion. In their model, the return to investment
is a function of the Brownian motion itself rather than its drift,
so the optimal policy is path-dependent. Lastly, there are real options
models in which the rate of production at each point in time is another
control variable in addition to the stopping time. These problems
are handled by the Hamilton-Jacobi-Bellman formulation of stochastic
control theory; see, for example, \citet{Knudsen1998}, \citet{Duckworth2000},
\citet{Alvarez2001b}, and \citet{Kumar2004}. 

In addition to the uncertainty in the profit stream, there is another
complicating but salient feature in our model: exit is possible after
investment. Among the papers that include this feature, \citet{McDonald1985}
study the valuation of a manufacturing firm facing a stochastic price
for its output product using option pricing techniques. In their model,
the product price is a geometric Brownian motion, and the firm can
shutdown and re-open its plant without cost at any point in time.
In contrast, \citet{Dixit1989} considers fixed cost of entry and
exit. In his model, the firm can enter and exit the industry as many
times as the firm wishes, and the profit stream is a geometric Brownian
motion. He shows that it is optimal to invest if the profit rate is
above an upper threshold and exit if it is below a lower threshold.
He performs a numerical comparative statics analysis and finds that
the upper (lower) threshold increases (decreases) in the volatility.
In his model, the investment (entry) decision can be exercised only
by an inactive firm; of course, the exit decision can be exercised
only by active firms. Our paper studies investment and exit decisions
in a quite different model: the firm has one opportunity to invest
in its operations while being active in the industry, and it can exit
at any point in time. Moreover, our comparative statics results are
analytical.

In the literatures on technology adoption and on exit, there is a
paucity of work on investment when the firm faces a declining profit
stream. To our knowledge, the current paper is the first to study
the impact of uncertainty on investment in an on-going project with
an exit option available both before and after an investment.

\section{The Basic Model\label{sec:basic-model}\label{sec:Optimal-exit}}

Consider a manufacturing firm whose product is produced with an aging
technology or process. Because of obsolescence,  its profit stream
is in decline (perhaps because a substitute product produced with
a new technology is encroaching upon the market). At any point in
time, the firm can stop the project by permanently closing its production
plant. 

The firm, seeking to maximize the expected discounted value of its
profit stream over an infinite horizon, must determine the best time
to cease operations and exit the market. The firm's profit rate at
time $t$ is a random variable $X_{t}$ where $\{X_{t}:t\ge0\}$ is
a stochastic process with continuous sample paths whose law of motion
we will specify shortly. We refer to $\{X_{t}:0\le t\le\tau\}$ as
the firm's \emph{profit stream} where the stopping time $\tau\le\infty$
is the time of exit. 

We model the firm's profit stream as a Brownian motion with constant
drift $\mu$ and volatility $\sigma$. Specifically, let $X_{t}$
denote the profit at time $t$ with $X_{t}=X_{0}+\nu t+\sigma B_{t}$
where $\{B_{t}:t\ge0\}$ is a one-dimensional standard Brownian motion,
so the profit stream has constant drift $\nu$ and constant volatility
$\sigma$. We pay particular attention to the case $\nu<0$ because
our main focus is modeling a declining profit stream. If the firm
begins operations at time $0$ and exits at a stopping time $\tau$,
the discounted value of its profit stream is $\int_{0}^{\tau}e^{-\alpha t}X_{t}dt$,
where $\alpha$ is the strictly positive discount rate. (To be more
precise, $\{X_{t}:t\ge0\}$ is a one-dimensional Brownian motion adapted
to a filtration $\{\mathcal{F}_{t}\}$ of a probability space $(\Omega,\mathcal{F},P)$.
The random variable $\tau$ is an element of $\mathcal{T}$, the set
of all non-negative stopping times with respect to the filtration
$\{\mathcal{F}_{t}\}$.)

To illustrate, suppose that the demand $D_{t}$ per unit time for
the firm's product is a Brownian motion with drift $\nu/p$, where
$p$ is the sales price per unit, and let $c$ be the fixed cost of
operation per unit time. Then the relationship between the demand
and the profit stream is linear: 
\begin{equation}
X_{t}=pD_{t}-c\:.\label{eq:X-D}
\end{equation}

Next, we review the solution and a few characteristics of the optimal
solution to the basic model. Our results below can be derived from
the more general results of \citet{Alvarez2001a} and \citet{Alvarez2003}.
In our basic model, at each point in time, the firm must elect either
to continue operations or to exit irrevocably. The firm seeks a stopping
time $\tau$ which maximizes 
\begin{equation}
E^{x}[\int_{0}^{\tau}X_{t}e^{-\alpha t}dt]\:,\label{eq:Exit-problem}
\end{equation}
where $E^{x}[\cdot]\equiv E[\cdot\vert X_{0}=x]$, the expectation
conditioned on $X_{0}=x$. The initial profit rate $x$ can be any
real number.

The firm receives cumulative discounted profit from continuation and,
implicitly, a reward of zero after exit. Hence, if the firm were to
exit at $\tau=0$, then it receives zero profit. For convenience,
we transform the problem to one with a reward from exit using the
strong Markov property of $X$ and the fact that $E^{x}\int_{0}^{\infty}\vert X_{t}\vert e^{-\alpha t}dt<\infty$
(\citealt{Alvarez2001a}, p. 318):
\begin{align}
E^{x}[\int_{0}^{\tau}X_{t}e^{-\alpha t}dt]= & E^{x}[\int_{0}^{\infty}X_{t}e^{-\alpha t}dt-\int_{\tau}^{\infty}X_{t}e^{-\alpha t}dt]\nonumber \\
= & E^{x}[\int_{0}^{\infty}X_{t}e^{-\alpha t}dt]-E^{x}[E^{X_{\tau}}[\int_{0}^{\infty}X_{t}e^{-\alpha t}dt]]\nonumber \\
= & \alpha^{-1}(x+\nu/\alpha)-\alpha^{-1}E^{x}[e^{-\alpha\tau}(X_{\tau}+\nu/\alpha)]\:,\label{eq:exit-problem-prime}
\end{align}
 which conforms to the standard stopping time problem with a reward
from exit. Here we used the identity
\[
E^{x}[\int_{0}^{\infty}X_{t}e^{-\alpha t}dt]=\alpha^{-1}(x+\nu/\alpha)\:,
\]
which can be derived using the Green function method (\citealt{Alvarez2001a},
p. 319). 

If the exit policy is stationary, then we can represent the stopping
time as $\tau_{D}$ which denotes the \emph{first exit time of the
process $X_{t}$ from a measurable set} $D\subset\mathbb{R}$:
\[
\tau_{D}\equiv\inf\{t>0:X_{t}\not\in D\}\:.
\]
In other words, the policy represented by the stopping time $\tau_{D}$
is to continue operations as long as $X_{t}\in D$ and stop when $X_{t}\not\in D$.
The set $D$ is called a \emph{continuation region}. The objective
function in Eq. (\ref{eq:Exit-problem}) or Eq. (\ref{eq:exit-problem-prime})
has no time-dependence other than through the process $X_{t}$ and
the discount factor $e^{-\alpha t}$; hence, we can show directly
(or use the argument of \citealt{Oksendal2003}, p. 220) that the
optimal policy, \emph{if it exists}, is stationary. Throughout this
section, we set $V(x;\nu)\equiv\sup_{\tau}E^{x}[\int_{0}^{\tau}X_{t}e^{-\alpha t}dt]$.
If $V(x;\nu)=E^{x}[\int_{0}^{\tau_{D}}X_{t}e^{-\alpha t}dt]$, then
$D$ is the \emph{optimal continuation region}, and an optimal policy
exists. 

\begin{proposition} \label{prop:exit-policy} The optimal policy
for the exit problem in Eq. (\ref{eq:Exit-problem}) exists, and the
optimal return is given by
\begin{equation}
V(x;\nu)=\left\{ \begin{array}{ll}
x/\alpha+\nu/\alpha^{2}+\frac{\sigma^{2}/\alpha}{\nu+\sqrt{\nu^{2}+2\alpha\sigma^{2}}}\exp[\frac{-\nu-\sqrt{\nu^{2}+2\alpha\sigma^{2}}}{\sigma^{2}}(x-\xi(\nu))] & \quad\textrm{if}\:x>\xi(\nu)\;,\\
0 & \quad\mbox{otherwise},
\end{array}\right.\label{eq:V0sol}
\end{equation}
where $\xi(\nu)$ is a negative number given by 
\begin{equation}
\xi(\nu)=-\frac{\nu}{\alpha}-\frac{\sigma^{2}}{\nu+\sqrt{\nu^{2}+2\alpha\sigma^{2}}}=\frac{\sigma^{2}}{\nu-\sqrt{\nu^{2}+2\alpha\sigma^{2}}}\:.\label{eq:xi0}
\end{equation}
Moreover, the optimal continuation region is $(\xi(\nu),\infty)$.
\end{proposition} 

\noindent \textbf{Proof}: The proof follows directly from Proposition
2 of \citet{Alvarez2001a}. In the terminology used by \citet{Alvarez2001a},
the decreasing fundamental solution is $e^{\phi x}$ where $\phi$
is defined in Eq. (\ref{eq:psi-phi}) below, the profit rate function
is $x$, and the reward from exit is zero. From the fact that $E[\int_{0}^{\infty}X_{t}e^{-\alpha t}dt]=x/\alpha+\nu/\alpha^{2}$,
it is straightforward to show that the optimal threshold is $\arg\max_{x}[-(x/\alpha+\nu/\alpha^{2})/e^{\phi x}]=\xi(\nu)$
and to verify that the optimal return function given by Proposition
2 of \citet{Alvarez2001a} reduces to Eq. (\ref{eq:V0sol}). \eproof

The most conventional way to solve a stopping time problem such as
the one in Eq. (\ref{eq:exit-problem-prime}) is to propose a candidate
function $V(\cdot;\nu)$ which is the return function from a candidate
policy and verify that the candidate function $V(\cdot;\nu)$ satisfies
a number of conditions as given by Theorem 10.4.1 of \citet{Oksendal2003}.
One of the conditions is that $V(\cdot;\nu)$ should satisfy a second-order
ordinary differential equation (ODE): $(-\alpha+\nu\partial_{x}+\frac{1}{2}\sigma^{2}\partial_{x}^{2})V(\cdot;\nu)=-x$.
Then $V(\cdot;\nu)$ is a sum of the term $(x/\alpha+\nu/\alpha^{2})$
and a function $f(x)$ which satisfies the ODE $(-\alpha+\nu\partial_{x}+\frac{1}{2}\sigma^{2}\partial_{x}^{2})f(x)=0$.
The general solution $f(x)$ is a linear combination of $e^{\psi x}$
and $e^{\phi x}$ where
\begin{equation}
\psi(\nu)=(-\nu+\sqrt{\nu^{2}+2\alpha\sigma^{2}})/\sigma^{2}\qquad\mbox{and}\qquad\phi(\nu)=(-\nu-\sqrt{\nu^{2}+2\alpha\sigma^{2}})/\sigma^{2}\;.\label{eq:psi-phi}
\end{equation}
 (We thank an anonymous referee for pointing us to Proposition 2 of
\citet{Alvarez2001a}: our exit problem in Eq. (\ref{eq:exit-problem-prime})
is a special case of Proposition 2 of \citet{Alvarez2001a} and admits
a much shorter proof above.) Interestingly, the optimal policy exists,
and the closed form solution $V(\cdot;\nu)$ is obtained. 

Proposition \ref{prop:exit-policy} says that the optimal policy is
to exit when the profit rate goes below the threshold $\xi(\nu)$.
It is intuitively clear that the firm will exit if its profit rate
has deteriorated below some threshold, but the fact that the threshold
is negative is not obvious. The reason $\xi(\nu)<0$ is that there
is value in waiting before taking an irrevocable action: even if $\nu<0$
and the current profit rate is slightly negative, it is possible for
the profit rate to turn positive in the future. If the profit stream
were deterministic and monotonically decreasing, then it would be
optimal to exit when the profit rate hits zero. This intuition regarding
the value of waiting is consistent with the fact that $\xi(\nu)$
increases to 0 as $\sigma\rightarrow0$, which follows from Eq. (\ref{eq:xi0})
when $\nu<0$. 

Now that we have a closed form solution for $V(\cdot;\cdot)$ and
$\xi(\cdot)$ in terms of all of the model parameters, their comparative
statics are straightforward as follows.

\begin{proposition} \label{prop:V-0-sigma} For any $\nu<0$, (i)
$\psi(\nu)$ decreases in $\sigma$ and $\nu$ and increases in $\alpha$;
$\phi(\nu)$ increases in $\sigma$ and decreases in $\nu$ and $\alpha$.
(ii) The exit threshold $\xi(\nu)$ decreases in $\sigma$ and $\nu$,
and it increases in $\alpha$. (iii) The optimal return $V(x;\nu)$
increases in $\sigma$ and $\nu$ for $x>\xi(\nu)$. \end{proposition}
(The proof is in the e-companion to this paper.)

The fact that $\xi(\nu)$ decreases in $\sigma$ is shown by \citet{Alvarez2003}
in a more general stopping time problem. The comparative statics of
$V(\cdot;\nu)$ with respect to $\sigma$ also follows from a more
general characteristic of stopping time problems obtained by Theorem
5 of \citet{Alvarez2003}: as $\sigma$ increases, there is more noise
in the profit stream, so there is a larger upturn potential as well
as a larger downturn risk. However, the firm can take advantage of
the upturn potential while avoiding downturn risk by exit. Hence,
the return function increases in $\sigma$. Similarly, because an
increase in $\nu$ improves the profit stream $X_{t}$, the return
function also increases for each continuation region $D$, so $V(\cdot;\nu)$
increases in $\nu$. 

In the next section, we consider a \emph{one-time} opportunity to
invest. The firm initially begins with drift $\mu<0$, and the investment
boosts the profit rate by $b$ and the drift by $\delta\ge0$. Because
there is only one investment opportunity, the post-investment problem
reduces to the basic model with drift $\nu=\mu+\delta$. For convenience,
we define
\[
V_{0}(\cdot)\equiv V(\cdot;\mu)\:
\]
as the optimal return to the exit problem with drift $\mu$,
\begin{equation}
\gamma_{p}\equiv\psi(\mu)=(-\mu+\sqrt{\mu^{2}+2\alpha\sigma^{2}})/\sigma^{2}\quad\mbox{and}\quad\gamma_{n}\equiv\phi(\mu)=(-\mu-\sqrt{\mu^{2}+2\alpha\sigma^{2}})/\sigma^{2}\;,\label{eq:gamma_pn}
\end{equation}
and 
\[
\xi_{0}\equiv\xi(\mu)
\]
as the optimal exit threshold with drift $\mu$. Similarly, we define
\[
V_{0}^{+}(\cdot)\equiv V(\cdot;\mu+\delta)\:
\]
as the optimal return to the exit problem with the boosted drift $\mu+\delta$,
and 
\begin{align}
\lambda\equiv\phi(\mu+\delta)= & [-(\mu+\delta)-\sqrt{(\mu+\delta)+2\alpha\sigma^{2}}]/\sigma^{2}\;,\label{eq:def-lambda}\\
\xi_{1}\equiv\xi(\mu+\delta)= & -\frac{\mu+\delta}{\alpha}+\frac{1}{\lambda}\:,\nonumber 
\end{align}
where $\xi_{1}$ is the post-investment optimal exit threshold. Since
the comparative statics obtained in Proposition \ref{prop:V-0-sigma}
applies for any negative drift, it applies to $V_{0}(\cdot)$, $\xi_{0}$,
$V^{+}(\cdot)$,  $\xi_{1}$, $\gamma_{p}$, $\gamma_{n}$, and $\lambda$
as long as $\mu<0$ and $\mu+\delta<0$. In the next section, $\max\{V_{0}^{+}(x+b)-k,0\}$
plays the role of the reward function from investment, where $k$
is the cost of investment. By Proposition \ref{prop:V-0-sigma}, the
reward from investment increases in the volatility,  a salient feature
of our model.

\section{The Model with One Investment Opportunity\label{sec:Optimal-investment} }

Consider the possibility of a once-in-a-lifetime investment. For instance,
manufacturers of 14-inch disk drives can, despite the writing on the
wall, improve the performance (recording capacity) of 14-inch drives
in order to immediately boost demand in the higher-end mainframe computer
market (\citealt{Christensen2000}, p. 19). Of course, eventual exit
is inevitable when $\mu<0$. 

For analytical tractability, our model allows only one investment
opportunity. As suggested by \citet{Fine1989}, in practice, the firm
might have multiple opportunities for gradual improvement in the technology/process.
The impact of multiple investment opportunities is beyond the scope
of this paper. 

\subsection{The Model}

We now include a one-time opportunity to implement an innovation which
improves the quality of the product or the process. The implementation
cost is $k>0$. If the quality of the product improves, then the demand
for the product increases; moreover, the demand declines more slowly.
Specifically, the investment boosts the current profit rate by $b$
and increases the drift by $\delta$. In terms of the example specified
in Eq. (\ref{eq:X-D}) where the profit rate is linear in the demand
$D_{t}$, investment induces an increase of $b$ in $pD_{t}$ (or,
equivalently, a decrease of $b$ in $c$) and an increase of $\delta$
in $p\cdot dD_{t}/dt$. If the firm invests at time $\tau$,  then
the improved profit stream follows the process 
\[
Y_{t}\equiv X_{t}+\delta(t-\tau)+b\:,\quad\mbox{for }\:t>\tau
\]
so that $dY_{t}=(\mu+\delta)dt+\sigma dB_{t}$. 

Prior to investment, the firm needs to find the optimal stopping time
$\tau$ at which to invest or to exit, whichever action results in
a better payoff. If the firm invests at time $\tau$, then its expected
return starting at time $\tau$ is $V_{0}^{+}(X_{\tau}+b)-k$ because
its expected cumulative profit stream after investment is $V_{0}^{+}(X_{\tau}+b)$
and the cost of investment is $k$. On the other hand, if the firm
exits at time $\tau$, then its return starting at time $\tau$ is
0. Hence, the firm receives the expected payoff of $\max\{V_{0}^{+}(X_{\tau}+b)-k,0\}$
at time $\tau$ when it makes its investment or exit decision. 

Let $x^{+}$ be the unique number which satisfies 
\begin{equation}
V_{0}^{+}(x^{+}+b)=k\quad.\label{eq:x-plus}
\end{equation}
(This definition uniquely determines $x^{+}$ because $V_{0}^{+}(x)$
is strictly increasing in $x$ for all $x$ such that $V_{0}^{+}(x)>0$.)
Then, at a stopping time $\tau$, it is optimal to exit if $X_{\tau}<x^{+}$
and invest if $X_{\tau}>x^{+}$ because $V_{0}^{+}(x+b)-k>0$ if $x>x^{+}$
and $V_{0}^{+}(x+b)-k<0$ if $x<x^{+}$. If the current profit rate
$X_{t}$ is $x^{+}$, then immediate investment and immediate exit
both yield zero expected return. Hence, our objective is to find
\begin{equation}
\bar{V}(x)\equiv\sup_{\tau\in\mathcal{T}}E^{x}[\int_{0}^{\tau}e^{-\alpha t}X_{t}dt+e^{-\alpha\tau}h(X_{\tau})]\:,\label{eq:supremum}
\end{equation}
where 
\begin{equation}
h(x)=\max\{0,V_{0}^{+}(x+b)-k\}\:\label{eq:htx}
\end{equation}
is the lump sum payoff when $x$ is the state when stopping occurs.

Next, we examine the conditions under which it is never optimal to
invest. Define
\begin{equation}
g\equiv\alpha(\int_{0}^{\infty}(b+\delta t)e^{-\alpha t}dt-k)=b+\delta/\alpha-k\alpha\label{eq:g-def}
\end{equation}
so that $g/\alpha$ is the net discounted gain from investment if
exit never occurs. 

\begin{proposition} \label{prop:set-I} Investment is never optimal
if and only if $g\le0$. \end{proposition}

\noindent \textbf{Proof}: To prove this proposition, we study the
difference between the return from immediate investment ($V_{0}^{+}(x+b)-k$)
and the optimal return from waiting to exit ($V_{0}(x)$) for $x>x^{+}$:
\begin{equation}
[V_{0}^{+}(x+b)-k]-V_{0}(x)=\frac{g}{\alpha}-\frac{1}{\lambda\alpha}\exp[\lambda(x+b-\xi_{1})]+\frac{1}{\gamma_{n}\alpha}\exp[\gamma_{n}(x-\xi_{0})]\:.\label{eq:Delta-V}
\end{equation}

\noindent Note that $-\frac{1}{\lambda\alpha}\exp[\lambda(x+b-\xi_{1})]<-\frac{1}{\gamma_{n}\alpha}\exp[\gamma_{n}(x-\xi_{0})]$
because $\lambda<\gamma_{n}<0$ and $\xi_{1}<\xi_{0}$ by Proposition
\ref{prop:V-0-sigma} (i) and (ii). If $g\le0$, then the right-hand-side
of Eq. (\ref{eq:Delta-V}) is negative for all $x>x^{+}$. Thus, at
any given stopping time $\tau$, the decision-maker is better off
choosing $V_{0}(X_{\tau})$ (return from waiting to exit) than $V_{0}^{+}(X_{\tau}+b)-k$
(return from immediate investment). We conclude that investment is
never optimal if $g\le0$. 

Now suppose that investment is also never optimal for $g>0$. In this
case, Eq. (\ref{eq:Delta-V}) is positive for sufficiently large $x$
because the two exponential terms converge to zero as $x\rightarrow\infty$.
Hence, at any $\tau$ such that $X_{\tau}$ is sufficiently large,
immediate investment is a better option than waiting to exit. This
contradicts the assumption that investment is never optimal. We conclude
that the policy of waiting to exit without ever investing is not optimal
when $g>0$. \eproof \\
In light of Proposition \ref{prop:set-I}, we assume $g>0$ for the
remainder of the paper unless otherwise noted. 

\subsection{Existence and Characterization of the Optimal Policy and the Optimal
Return  \label{subsec:Analysis}}

In this subsection, we verify that an optimal policy always exists
and that it is essentially unique. We show that, prior to investment,
the optimal policy is completely characterized by a pair $(\xi_{E},\xi_{I})$
of thresholds: exit if $X_{t}\le\xi_{E}$ and invest if $X_{t}\ge\xi_{I}$;
after investment, the problem reverts to the one analyzed at the end
of Sec. \ref{sec:Optimal-exit} where we established that it is optimal
to exit as soon as $X_{t}\le\xi_{1}$. Partial recompense for the
analytical difficulty implicit in our model is found in the closed
form solution for the optimal return function as given in Eq. (\ref{eq:V1-x}).

If the optimal policy exists, it is stationary because neither the
payoff $\max\{V_{0}^{+}(X_{t}+b)-k,0\}$ nor the profit stream has
any time-dependence other than through $X_{t}$ and $e^{-\alpha t}$.
Thus, it suffices to consider the class of stopping times $\tau_{D}=\inf\{t>0:X_{t}\not\in D\}$
expressed with respect to continuation regions $D$. Consequently,
we can express the objective function as 
\begin{equation}
R_{D}(x)=E^{x}[\int_{0}^{\tau_{D}}e^{-\alpha t}X_{t}dt+e^{-\alpha\tau_{D}}h(X_{\tau_{D}})]\;.\label{eq:V1tilde}
\end{equation}
In this new representation, the firm's policy is to continue operations
as long as $X_{t}\in D$ and to stop as soon as $X_{t}\not\in D$,
at which time the firm receives  $h(X_{t})$. Hence, if the optimal
policy exists, our objective is to find the optimal continuation region
$D^{*}$ such that 
\begin{equation}
\bar{V}(x)\equiv\sup_{D}R_{D}(x)=R_{D^{*}}(x)\;.\label{eq:V1}
\end{equation}

\begin{proposition}  \label{prop:g-positive} In the stopping time
problem of Eq. (\ref{eq:V1}), the optimal policy always exists; the
optimal expected return is uniquely given by
\begin{equation}
V_{1}(x)=\left\{ \begin{array}{ll}
x/\alpha+\mu/\alpha^{2}+a_{1}e^{\gamma_{p}x}+a_{2}e^{\gamma_{n}x}\,, & \;\mbox{for}\:x\in D^{*}=(\xi_{E},\xi_{I})\\
h(x)\quad, & \mbox{otherwise}
\end{array}\right.\:,\label{eq:V1-x}
\end{equation}
and the optimal continuation region is $D^{*}=(\xi_{E},\xi_{I})$:
it is optimal to exit when $x\le\xi_{E}$, invest when $x\ge\xi_{I}$,
and continue operations otherwise. \end{proposition}  The proof of
Proposition \ref{prop:g-positive} is in the e-companion to this paper.
The proof proceeds by considering $V_{1}(\cdot)$ as a candidate for
being the optimal return function and then verifying that $V_{1}(\cdot)$
satisfies all the sufficient conditions for being the optimal return
function specified in Theorem 10.4.1 of \citet{Oksendal2003}. (The
general necessary conditions for optimality of a return function with
a \emph{bounded} continuation region are obtained by \citealt{Alvarez2001a};
see also \citealt{Guo2001}.) Surprisingly, the existence proof is
rather difficult. It amounts to showing that the following boundary
conditions 
\begin{align}
V_{1}(\xi_{E}) & =\xi_{E}/\alpha+\mu/\alpha^{2}+a_{1}e^{\gamma_{p}\xi_{E}}+a_{2}e^{\gamma_{n}\xi_{E}}=h(\xi_{E})=0\;,\label{eq:bdry-1}\\
V_{1}(\xi_{I}) & =\xi_{I}/\alpha+\mu/\alpha^{2}+a_{1}e^{\gamma_{p}\xi_{I}}+a_{2}e^{\gamma_{n}\xi_{I}}\nonumber \\
 & =h(\xi_{I})=(\xi_{I}+b)/\alpha+\mu^{+}/\alpha^{2}-(\alpha\lambda)^{-1}e^{\lambda(\xi_{I}+b-\xi_{1})}-k\:,\label{eq:bdry-3}
\end{align}
along with the smooth-pasting conditions
\begin{align}
\partial_{x}V_{1}(\xi_{E}) & =\alpha^{-1}+\gamma_{p}a_{1}e^{\gamma_{p}\xi_{E}}+\gamma_{n}a_{2}e^{\gamma_{n}\xi_{E}}=\partial_{x}h(\xi_{E})=0\:,\label{eq:bdry-2}\\
\partial_{x}V_{1}(\xi_{I}) & =\alpha^{-1}+\gamma_{p}a_{1}e^{\gamma_{p}\xi_{I}}+\gamma_{n}a_{2}e^{\gamma_{n}\xi_{I}}=\partial_{x}h(\xi_{I})=\alpha^{-1}[1-e^{\lambda(\xi_{I}+b-\xi_{1})}]\:.\label{eq:bdry-4}
\end{align}
 always have a solution with desirable properties as stipulated by
Theorem 10.4.1 of \citet{Oksendal2003}. 

The optimal return function $V_{1}(\cdot)=\bar{V}(\cdot)$ is, of
course, unique. It follows that the optimal policy is unique: to stop
when $X_{t}\in\{x:x<\xi_{E}\:\mbox{or}\;x>\xi_{I}\}$ and continue
otherwise. If the current profit rate is $x^{+}$, then there is positive
probability that the profit rate will increase to a value bigger than
$x^{+}$ in the immediate future. Hence, the expected return from
waiting is positive, so $V_{1}(x^{+})>0$ and $x^{+}\in D^{*}=(\xi_{E},\xi_{I})$.
By Proposition \ref{prop:g-positive}, the firm's optimal policy is
to stop the first time $X_{t}$ hits $\xi_{E}$ or $\xi_{I}$ and
receive the reward $h(X_{t})$. Because $\xi_{E}<x^{+}<\xi_{I}$,
Eq. (\ref{eq:x-plus}) reveals that $V_{0}^{+}(\xi_{I}+b)-k>0$ and
$V_{0}^{+}(\xi_{E}+b)-k<0$. As anticipated, the firm's optimal action
at the stopping time $\tau_{D^{*}}$ depends on which end of the interval
$(\xi_{E},\xi_{I})$ the profit rate $X_{t}$ hits first. It is optimal
to exit if $X_{t}$ hits $\xi_{E}$ at time $\tau_{D^{*}}$, and it
is optimal to invest if $X_{t}$ hits $\xi_{I}$ at time $\tau_{D^{*}}$. 

\subsection{Comparative Statics\label{sec:Comparative-Statics}}

In this subsection, we effect a comparative statics analysis of $V_{1}(\cdot)$,
$\xi_{E}$, and $\xi_{I}$. We first establish the convexity of $V_{1}(\cdot)$,
which leads to the comparative statics with respect to $\sigma$. 

\begin{lemma} \label{lem:V1-convexity} The optimal return function
$V_{1}(\cdot)$ is convex. \end{lemma} (The proof is in the e-companion
to this paper.)

Next, we examine the comparative statics of $V_{1}(\cdot)$ with respect
to $\mu$ and $\sigma$. \begin{proposition}  \label{prop:V1-comp}
For all $x\in\mathbb{R}$, $V_{1}(x)$ is non-decreasing in $\mu$
and $\sigma$. In particular, $V_{1}(x)$ is strictly increasing in
$\mu$ for $x>\xi_{E}$. \end{proposition} 

\noindent \textbf{Proof:} To begin, note that $h(\cdot)$ is convex
and non-decreasing because $V_{0}^{+}(\cdot)$ is convex and non-decreasing.
Also note that $h(\cdot)$ is non-decreasing in both $\mu$ and $\sigma$
because $V_{0}(\cdot)$ is non-decreasing in $\mu$ and $\sigma$
by Proposition \ref{prop:V-0-sigma}. 

To show that $V_{1}(\cdot)$ is non-decreasing in $\mu$, let $V_{1}(x;\mu)$
and $h(x;\mu)$ denote the dependence of $V_{1}(x)$ and $h(x)$ on
the initial (pre-investment) drift $\mu$. Then for any $\beta>0$
and $x>\xi_{E}$, 
\begin{eqnarray*}
V_{1}(x;\mu) & = & E^{x}[\int_{0}^{T_{\mu}}X_{t}e^{-\alpha t}dt+e^{-\alpha T_{\mu}}h(X_{T_{\mu}};\mu)]\\
 & < & E^{x}[\int_{0}^{T_{\mu}}(X_{t}+\beta t)e^{-\alpha t}dt+e^{-\alpha T_{\mu}}h(X_{T_{\mu}}+\beta T_{\mu};\mu+\beta)]\le V_{1}(x;\mu+\beta)
\end{eqnarray*}
where $T_{\mu}$ is the optimal stopping time which maximizes $R_{D}(x)$
when the drift is $\mu$. In establishing the strict inequality, we
used the fact that $T_{\mu}>0$ for $x>\xi_{E}$, $h(x;\mu)$ is non-decreasing
in $x$ and $\mu$, and $T_{\mu}$ is suboptimal when the drift is
$\mu+\beta$. 

For comparative statics with respect to $\sigma$, we first assume
that $h(\cdot)$ does \emph{not} have functional dependence on $\sigma$.
By Lemma 2 and Theorem 4 of \citet{Alvarez2003} concerning the comparative
statics of more general stopping problems, because $V_{1}(\cdot)$
is convex and it is obtained as the return from stopping at $\tau_{(\xi_{E},\xi_{I})}$,
$V_{1}(\cdot)$ is non-decreasing in $\sigma$. (See also \citet{Ekstrom2004}).
Moreover, the reward function $h(\cdot)$ is non-decreasing in $\sigma$
by Prop \ref{prop:V-0-sigma}. Thus, $V_{1}(x)$ is non-decreasing
in $\sigma$\eproof  

By the same argument used in the proof of Theorem 5 of \citet{Alvarez2003},
the comparative statics of $\xi_{E}$ follows from Proposition \ref{prop:V1-comp}:

\begin{corollary} \label{cor:thresholds} The exit threshold $\xi_{E}$
satisfies $\partial_{\mu}\xi_{E}<0$ and $\partial_{\sigma^{2}}\xi_{E}\le0$.
\end{corollary}

\noindent \textbf{Proof:} Noting that $\xi_{E}=\inf\{x:V_{1}(x)>0\}$,
this result follows from the fact that $V_{1}(\cdot)$ is strictly
increasing in $\mu$ for $x>\xi_{E}$ and non-decreasing in $\sigma$
(by Proposition \ref{prop:V1-comp}). \hfill $\blacksquare$ 

In contrast, the comparative statics of $\xi_{I}$ is considerably
more complicated. Because $V_{1}(x)>V_{0}^{+}(x+b)-k$ if and only
if $x<\xi_{I}$, we can write $\xi_{I}=\sup\{x:V_{1}(x)-[V_{0}^{+}(x+b)-k]>0\}$.
Hence, the dependence of both $V_{1}(\cdot)$ and $V_{0}^{+}(\cdot)$
on $\mu$ and $\sigma$ determine the comparative statics of $\xi_{I}$.
This is in stark contrast to models in which the reward functions
have no $\sigma$ or $\mu$ dependence. (See, for example, Theorem
5 of \citealp{Alvarez2003}.) In order to examine the comparative
statics of $\xi_{I}$, we need to study the equations for both $\xi_{E}$
and $\xi_{I}$ as expressed by Eqs. (\ref{eq:y11}) and (\ref{eq:y12})
of Appendix B, where $\lambda$ is given by (\ref{eq:def-lambda}). 

Note that a closed-form expression for $\xi_{I}$ and $\xi_{E}$ can
not be obtained from Eqs. (\ref{eq:y11}) and (\ref{eq:y12}). In
contrast, using the closed-form expression for $\xi_{0}$, it was
straightforward to effect a complete comparative statics analysis
of $\xi_{0}$. Lack of a closed-form expression impairs our ability
to effect a comparative statics analysis of $\xi_{I}$. However, we
can obtain useful insights by examining the leading-order terms of
$\xi_{I}$ in power series expansions of $g$ when $b$ is close to
$\alpha k-\delta/\alpha$ ($g$ is small) and when $b$ is large ($g$
is large). We do not consider $\delta$ large because we restrict
our discussions to the interesting case $\mu^{+}<0$, i.e., the profit
stream is in decline even after investment. 

Using the expansions given in Lemmas \ref{lemma:small-g} and \ref{lemma:large-b}
of Appendix B, we obtain the limiting behavior of $\xi_{I}$ and $\xi_{E}$.
As $g\rightarrow0$, we find $\xi_{E}\rightarrow\xi_{0}$ and $\xi_{I}\rightarrow\infty$;
this echoes the intuition that it is almost never optimal to invest
when $g$ is close to zero. In the other limit where $b\rightarrow\infty$,
we find $\xi_{E}\rightarrow-\infty$ and $\xi_{I}-\xi_{E}\rightarrow0$;
this occurs because it is optimal to invest whenever $b$ is sufficiently
large. 

\begin{lemma} \label{lem:gamma-n} (i) $\partial_{\sigma^{2}}\gamma_{n}>0$
and $\partial_{\sigma^{2}}\lambda>0$; (ii) $\partial_{\mu}\gamma_{n}<0$
and $\partial_{\mu}\lambda<0$; (iii) $\lambda/\gamma_{n}\ge1$. \end{lemma} 

\noindent \textbf{Proof:} From Eqs. (\ref{eq:gamma_pn}) and (\ref{eq:def-lambda}),
(i) and (ii) can be shown after some algebra. Statement (iii) follows
from statement (ii) because $\gamma_{n}$ is equal to $\lambda$ if
$\mu$ is replaced by $\mu+\delta$ with $\delta\ge0$. \eproof

\begin{proposition}  \label{prop:comp-small-g} For sufficiently
small values of $g$, (i) $\partial_{\sigma^{2}}\xi_{I}>0$ and $\partial_{\sigma^{2}}\xi_{E}<0$,
and also (ii) $\partial_{\mu}\xi_{I}<0$. \end{proposition} 

\noindent \textbf{Proof:} Take partial derivatives of leading order
terms of Eqs. (\ref{eq:IIy10}) and (\ref{eq:IIy20}) with respect
to $\mu$ and $\sigma^{2}$ and use statement (i) of Proposition \ref{prop:V-0-sigma}.\hfill{}$\blacksquare$

\begin{proposition} \label{prop:comp-large-b} For $b$ sufficiently
large, (i) $\partial_{\sigma^{2}}\xi_{I}<0$ and $\partial_{\sigma^{2}}\xi_{E}<0$,
and also (ii) $\partial_{\mu}\xi_{I}<0$. \end{proposition}

\noindent \textbf{Proof:} (i) From the definition of $\xi_{0}$ and
Eq. (\ref{eq:y1-large-eta}), we have (a function $f(x)$ such that
$f(x)\rightarrow0$ as $x\rightarrow\infty$ is said to be $o(1)$)
\[
\partial_{\sigma^{2}}\xi_{E}=-\gamma_{n}^{-2}\partial_{\sigma^{2}}\gamma_{n}+\partial_{\sigma^{2}}\theta+o(1)=-z(e^{z}-1)^{-1}\lambda^{-2}\partial_{\sigma^{2}}\lambda+o(1)\:,
\]
where $\theta$ is a positive number defined by Eq. (\ref{eq:theta-eq.}),
$z\equiv-\lambda(\theta+\alpha k+\gamma_{n}^{-1}-\lambda^{-1})$ is
strictly positive by Lemma \ref{lem:theta-sign}, and $\partial_{\sigma^{2}}\theta$
is given by Eq. (\ref{eq:dtheta-dsigma}). Note that $z$ and $\theta$
are independent of $b$ so that they are not affected when we take
the limit as $b\rightarrow\infty$. Because $\partial_{\sigma^{2}}\lambda>0$
from Proposition \ref{prop:V-0-sigma} (i), we have $\partial_{\sigma^{2}}\xi_{E}<0$
for sufficiently large $b$. From Eq. (\ref{eq:y2-large-eta}), we
have $\partial_{\sigma^{2}}(\xi_{I}-\xi_{E})\rightarrow0$ as $b\rightarrow\infty$
so that 
\[
\partial_{\sigma^{2}}\xi_{I}=\partial_{\sigma^{2}}\xi_{E}+\partial_{\sigma^{2}}(\xi_{I}-\xi_{E})=-z(e^{z}-1)^{-1}\lambda^{-2}\partial_{\sigma^{2}}\lambda+o(1)\:.
\]
 Thus, $\partial_{\sigma^{2}}\xi_{I}<0$ for sufficiently large $b$.

(ii) By Corollary \ref{cor:thresholds}, $\partial_{\mu}\xi_{E}<0$.
From Eqs. (\ref{eq:y2-large-eta}) and (\ref{eq:dtheta-dmu}), 
\[
\partial_{\mu}\xi_{I}=\partial_{\mu}\xi_{0}+\partial_{\mu}(\xi_{E}-\xi_{0})+\partial_{\mu}(\xi_{I}-\xi_{E})=-\alpha^{-1}-z(e^{z}-1)^{-1}\lambda^{-2}\partial_{\mu}\lambda_{\mu}+o(1)\:.
\]
From the definition of $\lambda$ in Eq. (\ref{eq:def-lambda}) and
Proposition \ref{prop:V-0-sigma} (i), we have $-\alpha^{-1}-z(e^{z}-1)^{-1}\lambda^{-2}\partial_{\mu}\lambda_{\mu}<0$.
Thus, $\partial_{\mu}\xi_{I}<0$ for sufficiently large $b$. \hfill{}$\blacksquare$ 

In the conventional real options model, the reward function (corresponding
to $h(\cdot)$ in our problem) has no $\sigma$-dependence. In this
case, as shown by Theorems 6 and 7 of \citet{Alvarez2003},  the continuation
region is enlarged as $\sigma$ increases, so we anticipate that the
entry (exit) threshold increases (decreases) in the volatility. In
our model, however, $h(\cdot)$ has an explicit dependence on $\sigma$,
so the $\sigma$-dependence of the thresholds does not necessarily
follow the result by \citet{Alvarez2003}. When $g$ is small, $\partial_{\sigma^{2}}\xi_{E}<0$
and $\partial_{\sigma^{2}}\xi_{I}>0$: as the volatility $\sigma$
increases, it is optimal to wait longer to take advantage of the upturn
potential before taking an irreversible action. This is similar to
the result obtained numerically by \citet{Dixit1989} and proved analytically
by \citet{Alvarez2003}. However, when $b$ is large, Proposition
\ref{prop:comp-large-b} (i) asserts that $\partial_{\sigma^{2}}\xi_{E}<0$
and $\partial_{\sigma^{2}}\xi_{I}<0$. Notice that the result $\partial_{\sigma^{2}}\xi_{I}<0$
stands in contrast to the conventional intuition inherited from real
options theory. This counterintuitive result obtains because the return
from investment, $V_{0}^{+}(x+b)-k$, depends on $\sigma$. It is
important to note that the return from investment has dependence on
$\sigma$ only because exit is possible after investment. 

\begin{figure}[H]
\includegraphics[width=8.8cm]{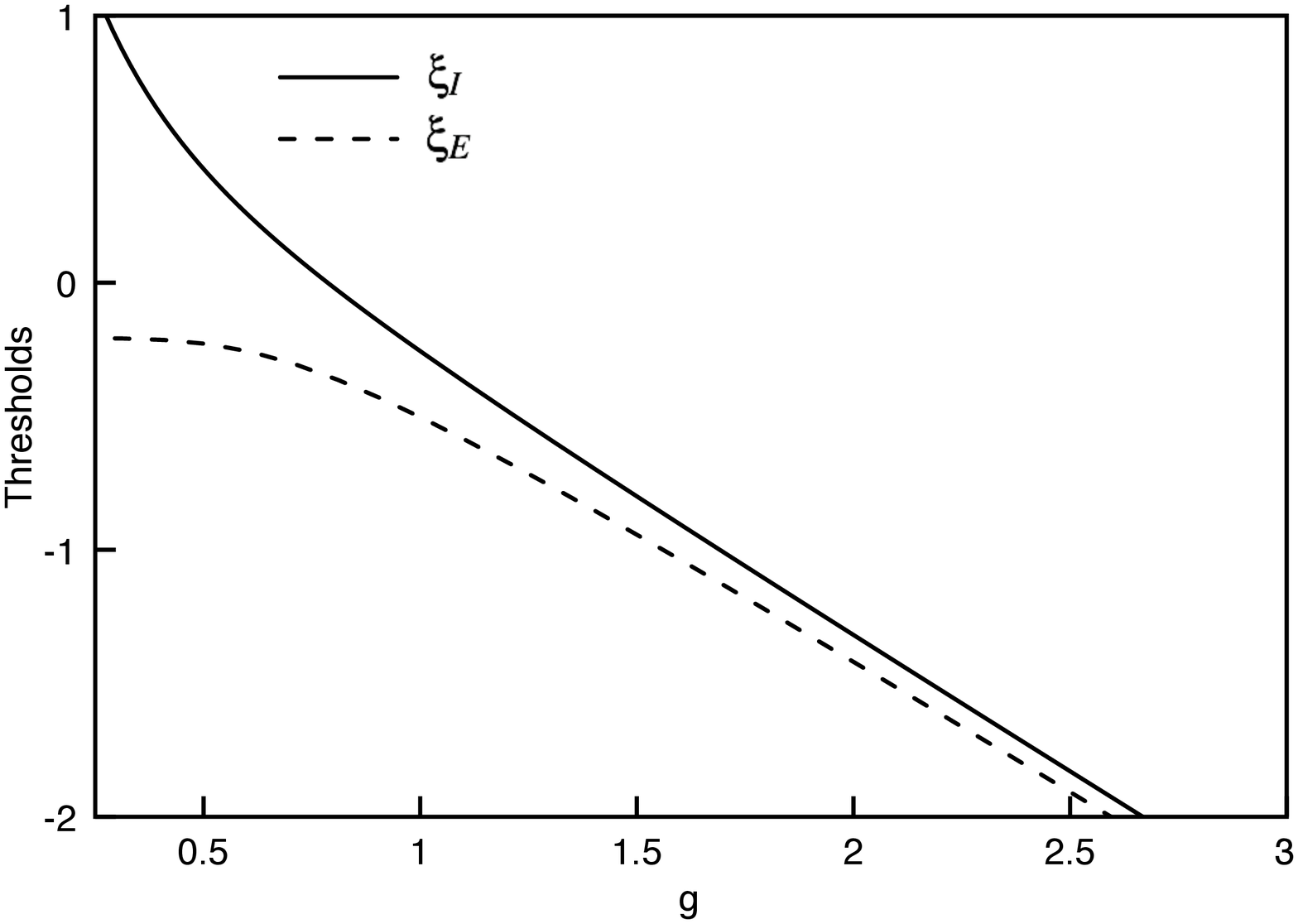}\caption{Thresholds as a function of $g$.}
\label{fig:Thresholds}
\end{figure}

\begin{figure}[H]
\includegraphics[width=8.8cm]{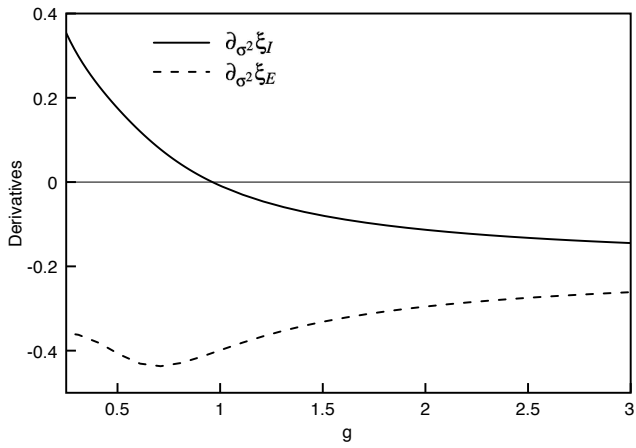}\caption{Derivatives of the thresholds with respect to $\sigma^{2}$ as a function
of $g$.}
\label{fig:Derivatives}
\end{figure}

The thresholds $\xi_{I}$ and $\xi_{E}$ and their comparative statics
with respect to $\sigma^{2}$ are illustrated by a numerical example
in Figs. \ref{fig:Thresholds} and \ref{fig:Derivatives}, in which
we set $\alpha=1,$ $\sigma^{2}=0.5$, $\mu=-1$, $\delta=0.1$, and
$k=0.5$. The graphs are shown as a function of $b+\delta/\alpha-\alpha k=b-0.4$.
 In Fig. \ref{fig:Derivatives}, notice that $\partial_{\sigma^{2}}\xi_{I}$
is positive for $g<0.96$ and negative for $g>0.96$. 

Another quantity of interest is the probability of investment prior
to the eventual exit and its dependence on volatility. Let $p_{I}(x)$
denote the probability (conditional on $X_{0}=x$ where $x\in(\xi_{E},\xi_{I})$)
that the profit rate hits $\xi_{I}$ before hitting $\xi_{E}$ (investment
is optimally made prior to exit). By II.4 and II.9 of \citet{Borodin2002},
\[
p_{I}(x)=\frac{\exp(-\frac{2\mu}{\sigma^{2}}x)-\exp(-\frac{2\mu}{\sigma^{2}}\xi_{E})}{\exp(-\frac{2\mu}{\sigma^{2}}\xi_{I})-\exp(-\frac{2\mu}{\sigma^{2}}\xi_{E})}=\frac{\exp[-\frac{2\mu}{\sigma^{2}}(x-\xi_{E})]-1}{\exp[-\frac{2\mu}{\sigma^{2}}(\xi_{I}-\xi_{E})]-1}\;.
\]

\begin{proposition} \label{prop:p-I}For sufficiently small and large
values of $g$, $p_{I}(\cdot)$ increases in $\sigma$. \end{proposition}

\noindent \textbf{Proof:} In the small-$g$ limit, by Eqs. (\ref{eq:IIy10})
and (\ref{eq:IIy20}), 
\[
p_{I}(x)=\exp[\frac{-2\mu}{\mu+\sqrt{\mu^{2}+2\alpha\sigma^{2}}}\log(g)]\cdot\{\exp[-\frac{2\mu}{\sigma^{2}}(x-\xi_{E})]-1\}\cdot(1+o(1))\;.
\]
Taking the derivative of the above with respect to $\sigma^{2}$,
we obtain
\[
\frac{d}{d\sigma^{2}}p_{I}(x)=p_{I}(x)\cdot\log(g)\frac{2\mu\alpha}{(\mu+\sqrt{\mu^{2}+2\alpha\sigma^{2}})^{2}\sqrt{\mu^{2}+2\alpha\sigma^{2}}}\cdot(1+o(1))\;.
\]
The leading-order term of $\frac{d}{d\sigma^{2}}p_{I}(x)$ is positive
because $\log(g)<0$ and $\mu<0$. 

Next, in the large-$g$ limit, by Eq. (\ref{eq:y2-large-eta}), both
$\xi_{I}-\xi_{E}$ and $x-\xi_{E}$ are bounded by $Cg^{-1}$ for
some positive constant $C$ because $x\in(\xi_{E},\xi_{I})$. Hence,
\[
p_{I}(x)=\frac{\exp[-\frac{2\mu}{\sigma^{2}}(x-\xi_{E})]-1}{\exp[-\frac{2\mu}{\sigma^{2}}(\xi_{I}-\xi_{E})]-1}=\frac{x-\xi_{E}}{\xi_{I}-\xi_{E}}(1+o(1))\;.
\]
 By Proposition \ref{prop:comp-large-b}, both $\xi_{I}$ and $\xi_{E}$
decrease in $\sigma$, and so
\[
\frac{d}{d\sigma^{2}}(\frac{x-\xi_{E}}{\xi_{I}-\xi_{E}})=-\partial_{\sigma^{2}}\xi_{I}\frac{x-\xi_{E}}{(\xi_{I}-\xi_{E})^{2}}-\partial_{\sigma^{2}}\xi_{E}\frac{\xi_{I}-x}{(\xi_{I}-\xi_{E})^{2}}>0\:.
\]
Hence, $\frac{d}{d\sigma^{2}}p_{I}(x)$ is positive for sufficiently
large $g$. \eproof

With a declining profit rate ($\mu<0$), the investment will only
be made if the profit rate is boosted by random noise. Hence, the
probability of making an investment before exit is expected to be
increasing in volatility. Although we suspect that this is a general
feature, the general comparative statics analysis is not available,
and we only can confirm the comparative statics in the two limiting
cases (small and large values of $g$) by Proposition \ref{prop:p-I}. 

\subsection{Extension: Switching to New Technology}

So far we have examined investment in the currently employed technology
which is becoming obsolete. Another important decision concerns when
to adopt the next-generation technology. In the context of the computer
hard-drive industry, some of the 14-inch drive manufacturers were
compelled to adopt 8-inch drives. In this subsection, we assume that
the decision-maker can adopt a new technology upon exit from the current
project and study how the exit value (expected profit from switching
to the new technology) affects the optimal policy. 

First, we examine the impact of adding a lump sum salvage value $s$
receivable at the time of exit. If plant and equipment are sold upon
exit, then we anticipate $s>0$. However, if there is employee severance
or liabilities associated with decommissioning of the business, then
$s<0$. 

\begin{lemma}  \label{lem:salvage}Let $V_{0}(\cdot;s)$ denote the
optimal return function when $s$ is the salvage value. Then
\[
V_{0}(x;s)=s+V_{0}(x-\alpha s)\:,
\]
and the exit threshold is $\xi(s)=\xi_{0}+\alpha s$. \end{lemma} 
(The proof is in the e-companion to this paper.)

Next, we consider a decision of whether and when to switch to a new
technology when there is a one-time investment opportunity to improve
the current technology. Without loss of generality, we assume that
there is no switching cost. Letting $s>0$ denote the expected cumulative
profit from switching to a new technology, we modify the investment
model considered in Sec. \ref{sec:Optimal-investment} by adding a
constant exit value $s$. Of course, even after investment in the
old technology, the firm can still switch to the new technology and
receive $s$. The objective is to find the optimal stopping time $\tau$
to maximize the following: 
\[
E^{x}[\int_{0}^{\tau}e^{-\alpha t}X_{t}dt+e^{-\alpha\tau}(h(X_{\tau}-\alpha s)+s)]\;.
\]

\begin{proposition}  \label{prop:new-tech} With the switching value
$s$, the optimal return from the investment decision problem is $V_{1}(x-\alpha s)+s$;
the investment threshold is $\xi_{I}+\alpha s$ and the switching
(exit) threshold is $\xi_{E}+\alpha s$. \end{proposition}  The proof
of Proposition \ref{prop:new-tech} is essentially the same as that
of Lemma \ref{lem:salvage}. With the opportunity to invest in a new
technology, the firm has less incentive to invest in or hold on to
the current technology. Hence, the thresholds for investment and switching
are higher when there is a profitable alternative technology. 

In this model of technology switching, we assumed that the cost of
switching is constant (or zero without loss of generality) and that
the expected profit from switching is independent of the current profit
rate $X_{t}$. In general, the switching cost and the expected profit
from switching may depend on the current profit rate if, for example,
$X_{t}$ represents the current demand (or market share) and the switched
technology serves the same market. There is also the possibility of
multiple and uncertain arrivals of improved technologies with different
dynamics of the profit stream. These complications are beyond the
scope of our paper. The problem of switching to future superior but
uncertain technologies have been studied by \citet{Alvarez2001} and
\citet{Balcer1984}, and the difficulty with switching to disruptive
technologies has been empirically studied by Christensen (1992 and
2000).

\section{Summary \label{sec:summary}}

Our analysis of investment under deteriorating conditions is congruent
with empirical reality as exemplified by the hard disk drive industry
and the many examples of obsolescent technologies enumerated in \citet{Christensen2000}
and \citet{Rosenberg1976}: it can be optimal to invest even in the
face of a declining profit stream and eventual displacement from the
market. Moreover, it can be optimal to remain in the market even if
the current profit rate is negative but above a threshold; it is optimal
to exit only when the profit rate has deteriorated sufficiently. 

In this paper, we studied a model of investment and exit decisions
under deteriorating conditions, and we proved that there exists the
optimal policy which is characterized by three thresholds: $\xi_{I},$
$\xi_{E}$, and $\xi_{1}$. Our comparative statics analysis with
respect to the volatility provided a novel and counterintuitive result.
As explained by \citet{Dixit1992}, illustrated by \citet{McDonald1986}
and \citet{Dixit1989}, and generalized by \citet{Alvarez2003}, it
is optimal to delay an irreversible action longer as the degree of
uncertainty increases in conventional real options models. In the
basic model of Sec. \ref{sec:Optimal-exit}, for instance, the exit
threshold $\xi_{0}$ always decreases in the volatility $\sigma$.
Similarly, in the model of Sec. \ref{sec:Optimal-investment}, the
exit threshold $\xi_{E}$ decreases in $\sigma$. The same intuition
suggests that $\xi_{I}$ increases in $\sigma$. Indeed, $\xi_{I}$
increases in $\sigma$ for sufficiently small $g$. However, we find
that $\xi_{I}$ decreases in $\sigma$ for sufficiently large $g$:
if the boost in the profit rate is sufficiently large, then it is
optimal to invest earlier as the uncertainty about the future profit
stream increases. (See Fig. \ref{fig:Derivatives}.) This counterintuitive
result is due to the firm's ability to control the time of its eventual
exit, a salient feature of our model. Because post-investment exit
is possible, the firm can take advantage of the volatility after investment,
so an increase in volatility induces an increase in the expected return
from investment and an increase in $\xi_{I}$ for sufficiently large
$g$.

\section*{Appendix: Equations for Thresholds}

Consider a solution $(\xi_{I},\xi_{E},a_{1},a_{2})$ to Eqs. (\ref{eq:bdry-1})
-- (\ref{eq:bdry-4}). For notational convenience, we define $\Delta_{IE}\equiv\xi_{I}-\xi_{E}$
and $\Delta_{E0}\equiv\xi_{E}-\xi_{0}\:.$ We eliminate $a_{1}$ and
$a_{2}$ from Eqs. (\ref{eq:bdry-1}) -- (\ref{eq:bdry-4}), and
we obtain 
\begin{eqnarray}
\Delta_{E0} & = & -ge^{-\gamma_{p}\Delta_{IE}}+(\lambda^{-1}-\gamma_{n}^{-1})e^{\lambda(\Delta_{IE}+\Delta_{E0}+b+\xi_{0}-\xi_{1})}e^{-\gamma_{p}\Delta_{IE}}\label{eq:y11}\\
 & = & -ge^{-\gamma_{n}\Delta_{IE}}+(\gamma_{p}^{-1}-\gamma_{n}^{-1})+(\lambda^{-1}-\gamma_{p}^{-1})e^{\lambda(\Delta_{IE}+\Delta_{E0}+b+\xi_{0}-\xi_{1})}e^{-\gamma_{n}\Delta_{IE}}\;,\label{eq:y12}
\end{eqnarray}
where $g$ is defined by Eq. (\ref{eq:g-def}).

In order to keep track of leading-order terms of power expansions
of $g$, we introduce a notation to denote the subleading order terms:
we say that $f(x)=o\left(j(x)\right)$ if $f(x)/j(x)\rightarrow0$
as $x\rightarrow0$, where $f(x)$ and $j(x)$ are functions of $x$.

\begin{lemma} \label{lemma:small-g} In the small-$g$ limit, 
\begin{eqnarray}
\Delta_{E0} & = & -g^{1-\gamma_{p}/\gamma_{n}}C(\delta)(1+o(1))\:,\label{eq:IIy10}\\
\Delta_{IE} & = & -\gamma_{n}^{-1}\ln(g^{-1})(1+o(1))\:,\label{eq:IIy20}
\end{eqnarray}
 where $C(\delta)=[(\gamma_{p}^{-1}-\gamma_{n}^{-1})]^{\gamma_{p}/\gamma_{n}}$
if $\delta>0$ and $C(\delta)=[(\gamma_{p}^{-1}-\gamma_{n}^{-1})(1-e^{\gamma_{n}b})]^{\gamma_{p}/\gamma_{n}}$
if $\delta=0$. \end{lemma} (The proof is  in the e-companion to
this paper.)

Similarly, we also say that $f(x)=o\left(j(x)\right)$ if $f(x)/j(x)\rightarrow0$
as $x\rightarrow\infty$.

\begin{lemma} \label{lemma:large-b} In the large-$b$ limit,
\begin{eqnarray}
\Delta_{E0} & = & -g+\theta+o(1)\label{eq:y1-large-eta}\\
\Delta_{IE} & = & -g^{-1}(\gamma_{p}\gamma_{n})^{-1}(1-\lambda\theta-\lambda/\gamma_{n})+o(g^{-1})\label{eq:y2-large-eta}
\end{eqnarray}
where $\theta$ is the unique positive solution to the equation
\begin{equation}
\theta=-\gamma_{n}^{-1}+\lambda^{-1}e^{\lambda(\theta+\alpha k-\delta/\alpha+\xi_{0}-\xi_{1})}\;.\label{eq:theta-eq.}
\end{equation}
\end{lemma} (The proof is in the e-companion to this paper.)

We need to obtain the comparative statics of $\theta$ in order to
examine the comparative statics of $\xi_{I}$ and $\xi_{E}$ in the
large-$b$ limit in Sec. \ref{sec:Comparative-Statics}. From Eq.
(\ref{eq:theta-eq.}) and the implicit function theorem, the partial
derivatives of $\theta$ with respect to $\sigma^{2}$ and $\mu$
are given by
\begin{align}
\partial_{\sigma^{2}}\theta= & \gamma_{n}^{-2}\partial_{\sigma^{2}}\gamma_{n}+\frac{e^{\lambda(\theta+\alpha k+\gamma_{n}^{-1}-\lambda^{-1})}}{1-e^{\lambda(\theta+\alpha k+\gamma_{n}^{-1}-\lambda^{-1})}}(\theta+\alpha k+\gamma_{n}^{-1}-\lambda^{-1})\lambda^{-1}\partial_{\sigma^{2}}\lambda\:,\label{eq:dtheta-dsigma}\\
\partial_{\mu}\theta= & \gamma_{n}^{-2}\partial_{\mu}\gamma_{n}+\frac{e^{\lambda(\theta+\alpha k+\gamma_{n}^{-1}-\lambda^{-1})}}{1-e^{\lambda(\theta+\alpha k+\gamma_{n}^{-1}-\lambda^{-1})}}(\theta+\alpha k+\gamma_{n}^{-1}-\lambda^{-1})\lambda^{-1}\partial_{\mu}\lambda\:.\label{eq:dtheta-dmu}
\end{align}
Finally, the following is a useful property of $\theta$:

\begin{lemma} \label{lem:theta-sign} $\theta+\alpha k+\gamma_{n}^{-1}-\lambda^{-1}>0$.\end{lemma}

\noindent \textbf{Proof}: For any value of $b$, $\Delta_{IE}+\Delta_{E0}+b+\xi_{0}-\xi_{1}>0$
always holds. In the limit $b\rightarrow\infty$, by Lemma \ref{lemma:large-b},
$\Delta_{IE}+\Delta_{E0}+b+\xi_{0}-\xi_{1}\longrightarrow\theta+\alpha k+\gamma_{n}^{-1}-\lambda^{-1}$,
which also must be positive. \eproof

\section*{Acknowledgments}

We thank Steven Lippman for many helpful suggestions and discussions.
We also thank two anonymous referees who provided suggestions which
considerably improved the manuscript. This research is partially supported
by the Price Center for Entrepreneurial Studies at UCLA and by the
UCLA Dissertation Year Fellowship.

\bibliographystyle{INFORMS2011}
\bibliography{OptInv1}

\vfill\pagebreak

{}

{}

\vfill\pagebreak

\section*{On-Line Appendix}

\section*{Appendix A: Proof of Proposition \ref{prop:g-positive}}

We first suppose that there is a solution ($\xi_{E},\xi_{I},a_{1},$
and $a_{2}$) to Eqs. (\ref{eq:bdry-1}) -- (\ref{eq:bdry-4}). (The
existence of a solution is assured by Lemma \ref{lem:existence}.)
Then Lemma \ref{lem:optimality} below establishes a condition under
which the optimal policy exists.

\begin{lemma}  \label{lem:optimality} Suppose that there exists
a solution ($\xi_{E},\xi_{I},a_{1},$ and $a_{2}$) to Eqs. (\ref{eq:bdry-1})
-- (\ref{eq:bdry-4}) that satisfies the constraints $\xi_{E}<x^{+}<\xi_{I}$
and
\begin{align}
x/\alpha+\mu/\alpha^{2}+a_{1}e^{\gamma_{p}x}+a_{2}e^{\gamma_{n}x} & \ge h(x)\quad\mbox{for}\quad x\in(\xi_{E},\xi_{I})\:,\label{eq:var-ineq}\\
\gamma_{p}^{2}a_{1}e^{\gamma_{p}x}+\gamma_{n}^{2}a_{2}e^{\gamma_{n}x} & \ge h^{\prime\prime}(x)\quad\mbox{for }\;x\in\{\xi_{E},\xi_{I}\}\:.\label{eq:convexity-ineq}
\end{align}
Then the optimal policy exists, and its expected return is given by
Eq. (\ref{eq:V1-x}); moreover, the optimal continuation region is
$D^{*}=(\xi_{E},\xi_{I})$.\end{lemma}  

\noindent \textbf{Proof}: Suppose there are $\xi_{E},\xi_{I},a_{1},$
and $a_{2}$ which satisfy Eqs. (\ref{eq:bdry-1})--(\ref{eq:bdry-4})
and Eq. (\ref{eq:var-ineq}). The infinitesimal generator for the
process $(t,X_{t})$ is given by $\partial_{t}+\mu\partial_{x}+\frac{1}{2}\sigma^{2}\partial_{x}^{2}$
(see \citealt{Oksendal2003}, p.222); however, when the time-dependence
of the return functions is only through the discount factor $e^{-\alpha t}$,
the infinitesimal generator can be conveniently replaced by 
\begin{equation}
\mathcal{A}\equiv-\alpha+\mu\partial_{x}+\frac{1}{2}\sigma^{2}\partial_{x}^{2}\;.\label{eq:A-operator}
\end{equation}

We consider $V_{1}(\cdot)$ defined in Eq. (\ref{eq:V1-x}) as a candidate
for the optimal return function. By Theorem 9.3.3, of \citet{Oksendal2003},
$V_{1}(x)$ is $R_{(\xi_{E},\xi_{I})}(x)$, the return function with
a continuation region $(\xi_{E},\xi_{I})$, because $V_{1}(x)$ satisfies
$\mathcal{A}V_{1}(x)=-x$ and the boundary conditions $V_{1}(\xi_{E})=h(\xi_{E})$
and $V_{1}(\xi_{I})=h(\xi_{I})$.

Next, we show that $V_{1}(\cdot)$ indeed coincides with $\bar{V}(\cdot)$
defined in Eq. (\ref{eq:supremum}) and establish that $(\xi_{E},\xi_{I})$
is the optimal continuation region. This is achieved by simply checking
the conditions of the variational inequalities given by Theorem 10.4.1
of \citet{Oksendal2003}.

First, as a preliminary condition for the variational inequalities,
$\{h^{-}(X_{\tau}):\tau\in\mathcal{T},\tau\le\infty\}$ is uniformly
integrable, which follows from the fact that $h^{-}(x)=0$. For a
second preliminary condition, we must show that $E^{x}[\int_{0}^{\infty}e^{-\alpha t}(X_{t})^{-}dt]<\infty$.
This follows from
\begin{eqnarray*}
E^{x}[\int_{0}^{\infty}e^{-\alpha t}(X_{t})^{-}dt]\le E^{x}[\int_{0}^{\infty}e^{-\alpha u}|X_{u}|du] & \leq & E^{x}[\int_{0}^{\infty}e^{-\alpha u}(|\mu u\vert+\vert x|+|\sigma B_{u}|)du]<\infty\:.
\end{eqnarray*}

Now we apply Theorem 10.4.1 of \citet{Oksendal2003} and find that
$V_{1}(x)\ge\bar{V}(x)$ because the following conditions are satisfied:
(a) $V_{1}(\cdot)$ is continuously differentiable {[}Eqs. (\ref{eq:bdry-2})
and (\ref{eq:bdry-4}){]} in $\mathbb{R}$, (b) $V_{1}(x)\ge h(x)$
for all $x\in\mathbb{R}$ {[}Eq. (\ref{eq:var-ineq}){]}, (c) $V_{1}(x)=h(x)$
for $x\in\{\xi_{E},\xi_{I}\}$ {[}Eqs. (\ref{eq:bdry-1}) and (\ref{eq:bdry-3}){]},
(d) $V_{1}(\cdot)$ is twice continuously differentiable except at
$\{\xi_{E},\xi_{I}\}$ (by the definition of $V_{1}(\cdot)$), (e)
the magnitudes of the second-order derivatives of $V_{1}(\cdot)$
are finite near $x=\xi_{E}$ and $\xi_{I}$, and (f) $\mathcal{A}V_{1}(x)\le-x$
for $x\in\mathbb{R}\backslash\{\xi_{E},\xi_{I}\}$. 

The last condition (f) has yet to be verified. By straightforward
algebra, $\mathcal{A}V_{1}(x)=-x$ for $x\in(\xi_{E},\xi_{I})$ and
\[
\mathcal{A}V_{1}(x)=\mathcal{A}h(x)=\left\{ \begin{array}{ll}
-x-g+\alpha^{-1}\delta e^{\lambda(x+b-\xi_{1})} & \textrm{if}\:x>\xi_{I}\;,\\
0 & \textrm{if}\:x<\xi_{E}
\end{array}\right.\;.
\]
By Eqs. (\ref{eq:convexity-ineq}), (\ref{eq:bdry-3}) and (\ref{eq:bdry-4}),
we have $\lim_{x\nearrow\xi_{I}}\mathcal{A}[V_{1}(x)-h(x)]\ge0$,
so $\lim_{x\nearrow\xi_{I}}\mathcal{A}V_{1}(x)=-\xi_{I}\ge\lim_{x\searrow\xi_{I}}\mathcal{A}h(x)=-\xi_{I}-g+\alpha^{-1}\delta e^{\lambda(\xi_{I}+b-\xi_{1})}\:.$
Since $e^{\lambda(x+b-\xi_{1})}$ decreases in $x$, we conclude that
$\mathcal{A}V_{1}(x)=\mathcal{A}h(x)<-x$ for $x>\xi_{I}$. Similarly,
by Eqs. (\ref{eq:convexity-ineq}), (\ref{eq:bdry-1}) and (\ref{eq:bdry-2}),
$\lim_{x\searrow\xi_{E}}\mathcal{A}[V_{1}(x)-h(x)]\ge0$ so that $-\xi_{E}\ge0$,
so $\mathcal{A}V_{1}(x)=\mathcal{A}h(x)=0<-x$ for $x<\xi_{E}$.

Finally, because $V_{1}(x)=R_{(\xi_{E},\xi_{I})}(x)\le\bar{V}(x)$,
we conclude that $V_{1}(x)=\bar{V}(x)$ is the unique optimal return
function and that $D^{*}=(\xi_{E},\xi_{I})$. \eproof

Next, we establish that there is a solution to Eqs. (\ref{eq:bdry-1})
-- (\ref{eq:bdry-4}) with some desirable properties. Note that a
solution to Eqs. (\ref{eq:y11}) and (\ref{eq:y12}) is also a solution
to Eqs. (\ref{eq:bdry-1}) -- (\ref{eq:bdry-4}). 

\begin{lemma} \label{lem:existence} There always exists a solution
($\xi_{E}$ and $\xi_{I}$) to Eqs. (\ref{eq:y11}) and (\ref{eq:y12})
that satisfies the constraints $\xi_{I}-\xi_{E}>0$ and $\xi_{I}+b-\xi_{1}>0$.
\end{lemma}

\noindent \textbf{Proof}: For fixed $\Delta_{IE}$, Eq. (\ref{eq:y11})
has a unique solution for $\Delta_{E0}$ because the left-hand-side
(LHS) is increasing while the right-hand-side (RHS) is decreasing
in $\Delta_{E0}$. Let $\Delta_{E0}=f(\Delta_{IE})$ be the solution
to Eq. (\ref{eq:y11}) as a function of $\Delta_{IE}$. 

At $\Delta_{IE}=0$, we have $f(0)+g=(\lambda^{-1}-\gamma_{n}^{-1})e^{\lambda(f(0)+b+\xi_{0}-\xi_{1})}$.
Because $f(0)+b+\xi_{0}-\xi_{1}=f(0)+g+\alpha k-\lambda^{-1}+\gamma_{n}^{-1}$,
if $f(0)+b+\xi_{0}-\xi_{1}\le0$, then $f(0)+g=(\lambda^{-1}-\gamma_{n}^{-1})e^{\lambda(f(0)+b+\xi_{0}-\xi_{1})}\ge\lambda^{-1}-\gamma_{n}^{-1}$,
which contradicts $f(0)+g+\alpha k-\lambda^{-1}+\gamma_{n}^{-1}\le0$.
Hence, $f(0)+b+\xi_{0}-\xi_{1}>0$. This also means that RHS of Eq.
(\ref{eq:y12}) is larger than RHS of Eq. (\ref{eq:y11}) if we set
$\Delta_{E0}=f(0)$ and $\Delta_{IE}=0$. For large positive values
of $\Delta_{IE}$, on the other hand, RHS of Eq. (\ref{eq:y12}) is
less than RHS of Eq. (\ref{eq:y11}) irrespective of the value of
$f(\Delta_{IE})$. Hence, there is a value of $\Delta_{IE}>0$ at
which RHS of Eq. (\ref{eq:y11}) equals RHS of Eq. (\ref{eq:y11}).
Thus, there is a solution to the simultaneous equations (\ref{eq:y11})
and (\ref{eq:y12}) with the constraint $\Delta_{IE}>0$.

Suppose that the solution satisfies $\xi_{I}+b-\xi_{1}\le0$. From
$\xi_{E}<\xi_{I}$, it follows that $\xi_{E}+b-\xi_{1}<0$. Hence,
\begin{equation}
\Delta_{E0}=\xi_{E}-\xi_{0}<-b+\xi_{1}-\xi_{0}=-g-\alpha k+(\lambda^{-1}-\gamma_{n}^{-1})\:.\label{eq:Del-E0-ineq}
\end{equation}
We also observe that $\Delta_{E0}<0$ because $b>0$ and $\xi_{1}<\xi_{0}$.
From Eq. (\ref{eq:y11}) and by the assumption that $e^{\lambda(\Delta_{IE}+\Delta_{E0}+b+\xi_{0}-\xi_{1})}\ge1$,
we have 
\begin{align*}
0>\Delta_{E0} & =-ge^{-\gamma_{p}\Delta_{IE}}+(\lambda^{-1}-\gamma_{n}^{-1})e^{\lambda(\Delta_{IE}+\Delta_{E0}+b+\xi_{0}-\xi_{1})}e^{-\gamma_{p}\Delta_{IE}}\\
 & \ge e^{-\gamma_{p}\Delta_{IE}}[-g+(\lambda^{-1}-\gamma_{n}^{-1})]>-g+(\lambda^{-1}-\gamma_{n}^{-1})\:,
\end{align*}
where the last inequality holds because $\lambda^{-1}-\gamma_{n}^{-1}\ge0$
{[}by the fact that $\phi(\nu)$ decreases in $\nu$ from Proposition
\ref{prop:V-0-sigma} (i){]} and $e^{-\gamma_{p}\Delta_{IE}}<1$.
The inequality $\Delta_{E0}>-g+(\lambda^{-1}-\gamma_{n}^{-1})$ contradicts
Eq. (\ref{eq:Del-E0-ineq}). \eproof

Lemma \ref{lem:existence} ensures that there exists a solution to
Eqs. (\ref{eq:bdry-1}) -- (\ref{eq:bdry-4}) with $\xi_{I}>\xi_{E}$,
and the solution satisfies the inequality $e^{\lambda(\Delta_{IE}+\Delta_{E0}+b+\xi_{0}-\xi_{1})}<1$.
We have yet to prove that $\xi_{I}>x^{+}>\xi_{E}$ and that Eqs. (\ref{eq:var-ineq})
and (\ref{eq:convexity-ineq}) hold. For the rest of this Appendix,
we assume a solution ($\xi_{I}$, $\xi_{E}$, $a_{1}$, $a_{2}$)
to Eqs. (\ref{eq:bdry-1}) -- (\ref{eq:bdry-4}) that satisfies $\xi_{I}>\xi_{E}$
and $\xi_{I}+b-\xi_{1}>0$, and we consider a test function $\phi(\cdot)$
defined as follows: 
\[
\phi(x)=x/\alpha+\mu/\alpha^{2}+a_{1}e^{\gamma_{p}x}+a_{2}e^{\gamma_{n}x}\:.
\]

In the following Lemma, we establish a very convenient property of
exponential functions which will be used in the proofs of forthcoming
Lemmas. Recall that $\lambda$ is defined in Eq. (\ref{eq:def-lambda})
and that $\lambda<\gamma_{n}$ by Proposition \ref{prop:V-0-sigma}
(i).

\begin{lemma} \label{lem:exponential} Given any function of the
form
\[
f(x)=C_{1}\exp(\gamma_{p}x)+C_{2}\exp(\gamma_{n}x)+C_{3}\exp(\lambda x)
\]
where $C_{1}>0$, suppose (i) $C_{2}>0$ or (ii) $C_{2}<0$ and $C_{3}<0$.
If $f(y)>0$ for some $y$, then $f(x)>0$ for all $x>y$.\end{lemma} 

\noindent \textbf{Proof}: (i) If $C_{2}>0$ and $C_{3}>0$, then the
lemma is obvious. If $C_{2}>0$ and $C_{3}<0$, then $\vert C_{3}\exp(\lambda x)\vert$
is decreasing in $x$ at a faster rate than $C_{2}\exp(\gamma_{n}x)$
while $C_{1}\exp(\gamma_{p}x)$ is increasing in $x$, so the lemma
follows. (ii) If $C_{2}<0$ and $C_{3}<0$, then $\vert C_{2}\exp(\gamma_{n}x)+C_{3}\exp(\lambda x)\vert$
decreases in $x$ while while $C_{1}\exp(\gamma_{p}x)$ increases
in $x$, so the lemma follows again. \eproof 

\begin{lemma} \label{lem:a-positive} The coefficients $a_{1}$ and
$a_{2}$ are positive. \end{lemma} 

\noindent \textbf{Proof}: Given a solution to Eqs. (\ref{eq:bdry-1})
-- (\ref{eq:bdry-4}) with the conditions $\xi_{I}>\xi_{E}$ and
$\xi_{I}+b-\xi_{1}>0$, $a_{1}$ and $a_{2}$ can be obtained from
Eqs. (\ref{eq:bdry-3}) and (\ref{eq:bdry-4}):
\begin{align*}
a_{1} & =\alpha^{-1}\gamma_{p}^{-1}\frac{e^{\gamma_{n}\xi_{I}}}{e^{\gamma_{p}\xi_{I}+\gamma_{n}\xi_{E}}-e^{\gamma_{p}\xi_{E}+\gamma_{n}\xi_{I}}}(1-e^{\lambda(\xi_{I}+b-\xi_{1})}e^{-\gamma_{n}(\xi_{I}-\xi_{E})})\:,\\
a_{2} & =-\alpha^{-1}\gamma_{n}^{-1}\frac{e^{\gamma_{p}\xi_{I}}}{e^{\gamma_{p}\xi_{I}+\gamma_{n}\xi_{E}}-e^{\gamma_{p}\xi_{E}+\gamma_{n}\xi_{I}}}(1-e^{\lambda(\xi_{I}+b-\xi_{1})}e^{-\gamma_{p}(\xi_{I}-\xi_{E})})\:.
\end{align*}
The denominator $e^{\gamma_{p}\xi_{I}+\gamma_{n}\xi_{E}}-e^{\gamma_{p}\xi_{E}+\gamma_{n}\xi_{I}}$
is always positive because $\xi_{I}>\xi_{E}$. Moreover, because $e^{\lambda(\xi_{I}+b-\xi_{1})}<1$
and $e^{-\gamma_{p}(\xi_{I}-\xi_{E})}<1$, $a_{2}$ is strictly positive. 

Suppose $a_{1}\le0$. (i) If $\xi_{E}<\xi_{I}\le x^{+}$, then $\phi(\xi_{E})=0$
and $\phi(\xi_{I})\le0$, so the first-derivative $\phi^{\prime}(y)$
takes a negative value at some point $y$ in the interval $(\xi_{E},\xi_{I})$.
We also know that $\phi^{\prime}(\xi_{E})=0$, so the second derivative
$\phi^{\prime\prime}(x)$ takes a negative value somewhere in the
interval $(\xi_{E},y)$; this is only possible if $a_{1}<0$. By Lemma
\ref{lem:exponential}, $\phi^{\prime\prime}(x)$ takes a negative
value in the interval $(y,\xi_{I})$ so that $\phi^{\prime}(\xi_{I})<0$.
This contradicts the condition $\phi^{\prime}(\xi_{I})=\alpha^{-1}[1-e^{\lambda(\xi_{I}+b-\xi_{1})}]>0$.

(ii) Suppose $x^{+}\le\xi_{E}<\xi_{I}$, and consider the function
$f(x)\equiv\phi(x)-h(x)$. By Eqs. (\ref{eq:bdry-3}) and (\ref{eq:bdry-4}),
$f(\xi_{I})=0$ and $f^{\prime}(\xi_{I})=0$. Because $\phi^{\prime}(\xi_{E})=0$
and $h^{\prime}(\xi_{E})>0$, we have $f^{\prime}(\xi_{E})<0$. From
the functional form $f^{\prime}(x)=\gamma_{p}a_{1}e^{\gamma_{p}x}+\gamma_{n}a_{2}e^{\gamma_{n}x}+\alpha e^{\lambda(x+b-\xi_{1})}$
and by Lemma \ref{lem:exponential}, $f^{\prime}(x)<0$ for all $x>\xi_{E}$
because $\gamma_{p}a_{1}\le0$ and $\gamma_{n}a_{2}<0$, so $f(\cdot)$
strictly decreases for $x>\xi_{E}$. Hence $f(\xi_{I})=0$ is impossible. 

(iii) The only remaining case is $\xi_{E}<x^{+}<\xi_{I}$. As was
argued in case (i), if $\phi(x^{+})\le0$, then $\phi(x)$ decreases
in $x$ for all $x>x^{+}$, in which case $f(\xi_{I})=0$ cannot be
achieved. Suppose $\phi(x^{+})>0$. Then $f(x^{+})>0$ and $f(\xi_{I})=0$,
so there is some $y\in(x^{+},\xi_{I})$ such that $f^{\prime}(y)<0$.
By Lemma \ref{lem:exponential}, $f^{\prime}(x)<0$ for all $x>y$,
in which case $f^{\prime}(\xi_{I})=0$ is impossible. 

Because $a_{1}\le0$ is impossible in all possible cases (i)--(iii),
we conclude that $a_{1}>0$.\eproof

The following Lemma ensures that the solution $\xi_{E}$ and $\xi_{I}$
satisfies $\xi_{E}<x^{+}<\xi_{I}$, one of the constraints required
by Lemma \ref{lem:optimality}.

\begin{lemma} \label{lem:xi-order} The inequality $\xi_{E}<x^{+}<\xi_{I}$
is satisfied. \end{lemma} 

\noindent \textbf{Proof}: Define $f(x)\equiv\phi(x)-h(x)$. By Lemma
\ref{lem:a-positive}, $\phi(\cdot)$ is strictly convex and strictly
increasing for $x>\xi_{E}$ because $\phi^{\prime}(\xi_{E})=0$. Also,
$\phi(x)$ is  positive for $x>\xi_{E}$ because $\phi(\xi_{E})=0$. 

(i) Suppose $\xi_{E}<\xi_{I}\le x^{+}$. Then $\phi(x)>0$ for all
$x>\xi_{E}$, but $(\xi_{I}+b)/\alpha+\mu^{+}/\alpha^{2}-(\alpha\lambda)^{-1}e^{\lambda(\xi_{I}+b-\xi_{1})}-k\le0$
because $\xi_{I}\le x^{+}$. Hence, Eq. (\ref{eq:bdry-3}) cannot
be satisfied.

(ii) Suppose $x^{+}\le\xi_{E}<\xi_{I}$. Because $f(\xi_{E})=0$ and
$f(\xi_{I})>0$, there is some $y\in(\xi_{E},\xi_{I})$ such that
$f^{\prime}(y)>0$. By Lemmas \ref{lem:exponential} and \ref{lem:a-positive},
$f^{\prime}(x)>0$ for all $x>y$, which contradicts the condition
$f^{\prime}(\xi_{I})=0$.\eproof 

Finally, it remains to show Eqs. (\ref{eq:var-ineq}) and (\ref{eq:convexity-ineq}).

\begin{lemma} \label{lem:var-ineq}The constraints Eqs. (\ref{eq:var-ineq})
and (\ref{eq:convexity-ineq}) are satisfied.\end{lemma} 

\noindent \textbf{Proof}: Because $\phi(x)$ is positive for $x>\xi_{E}$,
we have $\phi(x)>h(x)=0$ for $x\in(\xi_{E},x^{+}]$. Now consider
the function $f(x)\equiv\phi(x)-h(x)$ for $x\in(x^{+},\xi_{I})$.
Suppose that $f(y)<0$ for some $y\in(x^{+},\xi_{I})$. Then $f^{\prime}(z)>0$
for some $z\in(y,\xi_{I})$ because $f(\xi_{I})=0$. By Lemma \ref{lem:exponential},
$f^{\prime}(x)>0$ for all $x>z$, which contradicts the condition
$f^{\prime}(\xi_{I})=0$. Therefore, $\phi(x)\ge h(x)$ for all $x\in(\xi_{E},\xi_{I})$.

\noindent Next, $\phi^{\prime\prime}(\xi_{E})>0=h^{\prime\prime}(\xi_{E})$
because $a_{1}>0$ and $a_{2}>0$. In the interval $(x^{+},\xi_{I})$,
$f(\cdot)$ decreased from $f(x^{+})>0$ to $f(\xi_{I})$ so there
is some $y\in(x^{+},\xi_{I})$ at which $f^{\prime}(y)<0$. Because
$f^{\prime}(\xi_{I})=0$, there is some $z\in(y,\xi_{I})$ at which
$f^{\prime\prime}(z)>0$. By Lemma \ref{lem:exponential}, $f^{\prime\prime}(\xi_{I})>0$.
Hence, Eq. (\ref{eq:convexity-ineq}) is satisfied. \eproof 

By Lemmas \ref{lem:optimality}, \ref{lem:existence}, \ref{lem:xi-order},
and \ref{lem:var-ineq}, Proposition \ref{prop:g-positive} is proved.

\section*{Appendix B}

\noindent \textbf{Proof of Proposition \ref{prop:V-0-sigma}:} (i)
From Eq. (\ref{eq:psi-phi}), it is straightforward to obtain the
following inequalities:
\begin{align*}
\frac{\partial\psi(\nu)}{\partial(\sigma^{2})}=\frac{-\nu^{2}-\alpha\sigma^{2}+\nu\sqrt{\nu^{2}+2\alpha\sigma^{2}}}{\sigma^{4}\sqrt{\nu^{2}+2\alpha\sigma^{2}}}<0\:, & \frac{\partial\phi(\nu)}{\partial(\sigma^{2})}=\frac{\nu^{2}+\alpha\sigma^{2}+\nu\sqrt{\nu^{2}+2\alpha\sigma^{2}}}{\sigma^{4}\sqrt{\nu^{2}+2\alpha\sigma^{2}}}>0\:,\\
\frac{\partial\psi(\nu)}{\partial\nu}=\frac{\nu-\sqrt{\nu^{2}+2\alpha\sigma^{2}}}{\sigma^{2}\sqrt{\nu^{2}+2\alpha\sigma^{2}}}<0\;, & \frac{\partial\phi(\nu)}{\partial\nu}=\frac{-\nu-\sqrt{\nu^{2}+2\alpha\sigma^{2}}}{\sigma^{2}\sqrt{\nu^{2}+2\alpha\sigma^{2}}}<0\\
\frac{\partial\psi(\nu)}{\partial\alpha}=\frac{1}{\sqrt{\nu^{2}+2\alpha\sigma^{2}}}>0\;, & \frac{\partial\phi(\nu)}{\partial\alpha}=-\frac{1}{\sqrt{\nu^{2}+2\alpha\sigma^{2}}}<0\qquad.
\end{align*}

(ii) We first recognize that $\xi(\nu)=-1/\psi(\nu)$. Using the chain
rule of differentiation ($\frac{\partial\xi(\nu)}{\partial z}=\frac{1}{\psi^{2}(\nu)}\frac{\partial\psi(\nu)}{\partial z}$
for any model parameter $z$), we can obtain the comparative statics
of $\xi(\nu)$ from the comparative statics of $\psi(\nu)$ above. 

(iii) For $x>\xi(\nu)$, we can express
\[
V(x;\nu)=\frac{x}{\alpha}+\frac{\nu}{\alpha^{2}}-\frac{1}{\alpha\phi(\nu)}\exp[\phi(\nu)(x-\xi(\nu))]\:.
\]
Hence, 
\begin{align*}
\frac{\partial V(x;\nu)}{\partial\sigma^{2}}= & \frac{1}{\alpha\phi^{2}(\nu)}\exp[\phi(\nu)(x-\xi(\nu))]\\
 & \times[\frac{\partial\phi(\nu)}{\partial\sigma^{2}}-\phi(\nu)\frac{\partial\phi(\nu)}{\partial\sigma^{2}}(x-\xi(\nu))-\frac{\phi^{2}(\nu)}{\psi^{2}(\nu)}\cdot\frac{\partial\psi(\nu)}{\partial\sigma^{2}}]\:.
\end{align*}
From the fact that $\phi(\nu)<0$, $x>\xi(\nu)$, $\frac{\partial\psi(\nu)}{\partial\sigma^{2}}<0$,
and 
\[
\frac{\partial\phi(\nu)}{\partial\sigma^{2}}=\frac{\nu^{2}+\alpha\sigma^{2}+\nu\sqrt{\nu^{2}+2\alpha\sigma^{2}}}{\sigma^{4}\sqrt{\nu^{2}+2\alpha\sigma^{2}}}>0\:,
\]
we conclude that $\frac{\partial V(x;\nu)}{\partial\sigma^{2}}>0$.

The comparative statics $\frac{\partial V(x;\nu)}{\partial\nu}>0$
is easier to see from the objective function in Eq. (\ref{eq:Exit-problem})
which increases in $\nu$ with a fixed $\tau$. Alternatively, $\frac{\partial V(x;\nu)}{\partial\nu}>0$
can be directly shown by algebra. \eproof

\noindent \textbf{Proof of Lemma \ref{lem:V1-convexity}:} Let $x_{\theta}=\theta x_{1}+(1-\theta)x_{2}$
where $\theta\in(0,1)$ and $x_{1}\neq x_{2}$, and let $\tau^{*}$
denote the optimal stopping time conditional on $X_{0}=x_{\theta}$.
From the fact that $h(\cdot)$ is convex and that 
\[
x_{\theta}+\mu t+\sigma B_{t}=\theta(x_{1}+\mu t+\sigma B_{t})+(1-\theta)(x_{2}+\mu t+\sigma B_{t})\:,
\]
we obtain the following inequality:
\begin{align*}
V_{1}(x_{\theta})= & E[\int_{0}^{\tau^{*}}e^{-\alpha t}(x_{\theta}+\mu t+\sigma B_{t})dt+e^{-\alpha\tau^{*}}h(x_{\theta}+\mu\tau^{*}+\sigma B_{\tau^{*}})]\\
\le & \theta E^{x_{1}}[\int_{0}^{\tau^{*}}e^{-\alpha t}X_{t}dt+e^{-\alpha\tau^{*}}h(X_{\tau^{*}})]\\
 & +(1-\theta)E^{x_{1}}[\int_{0}^{\tau^{*}}e^{-\alpha t}X_{t}dt+e^{-\alpha\tau^{*}}h(X_{\tau^{*}})]\\
\le & \theta V_{1}(x_{1})+(1-\theta)V_{1}(x_{2})\:.
\end{align*}
\eproof

\noindent \textbf{Proof of Lemma \ref{lem:salvage}:} Because $se^{-\alpha\tau}=s-\int_{0}^{\tau}\alpha se^{-\alpha t}dt$,
\begin{align*}
V_{0}(x;s)= & E^{x}[\int_{0}^{\tau}X_{t}e^{-\alpha t}dt+se^{-\alpha\tau}]E^{x}=s+E^{x}[\int_{0}^{\tau}(X_{t}-\alpha s)e^{-\alpha t}dt]\\
= & s+E^{x-\alpha s}[\int_{0}^{\tau}X_{t}e^{-\alpha t}dt]=s+V_{0}(x-\alpha s)\;.
\end{align*}

\noindent Because $V_{0}(x;s)$ is increasing in $s$, the exit threshold
is 
\[
\xi(s)=\inf\{x:V_{0}(x;s)>s\}=\xi_{0}+\alpha s\:.
\]
 \hfill$\blacksquare$

\noindent \textbf{Proof of Lemma \ref{lemma:small-g}:} First, we
notice that if $g=0$, then $\Delta_{E0}=\xi_{E}-\xi_{0}=0$ and $\Delta_{IE}=\xi_{I}-\xi_{E}=\infty$.
Hence, $\Delta_{E0}\rightarrow0$ and $\Delta_{IE}\rightarrow\infty$
as $g\rightarrow0$. 

Suppose that $\delta>0$. The first term of the right-hand-side (RHS)
of Eq. (\ref{eq:y11}) strictly dominates the second term so that
\begin{equation}
g^{-1}\exp[\lambda(\Delta_{IE}+\Delta_{E0}+b+\xi_{0}-\xi_{1})]\rightarrow0\quad\textrm{as}\;g\rightarrow0\;.\label{eq:domination1}
\end{equation}
(Otherwise, if the second term of RHS of Eq. (\ref{eq:y11}) is dominant,
then $\Delta_{E0}>0$ in the $g\rightarrow0$ limit; if both terms
converge to zero at the same rate as $g\rightarrow0$, then the RHS
of Eq. (\ref{eq:y12}) converges to $(\gamma_{p}^{-1}-\gamma_{n}^{-1})$
as $g\rightarrow0$.) From Eq. (\ref{eq:domination1}), the leading
order terms in RHS of Eq. (\ref{eq:y12}) are contained in the first
two terms: $-ge^{-\gamma_{n}\Delta_{IE}}+(\gamma_{p}^{-1}-\gamma_{n}^{-1})$
in agreement with Eq. (\ref{eq:IIy20}). From the fact that $\lim_{g\rightarrow0}\Delta_{E0}=0$,
the only possible leading order term of $\Delta_{IE}$ is $\gamma_{n}^{-1}\ln[g(\gamma_{p}^{-1}-\gamma_{n}^{-1})^{-1}]$.
The leading-order terms of $\Delta_{IE}=\gamma_{n}^{-1}\ln[g(\gamma_{p}^{-1}-\gamma_{n}^{-1})^{-1}]+o(1)$
is consistent with the condition in Eq. (\ref{eq:domination1}) because
$\lambda/\gamma_{n}>1$ by Proposition \ref{prop:V-0-sigma} (i).
Finally, using the leading-order term of $\Delta_{IE}$ in Eq. (\ref{eq:y11}),
we obtain Eq. (\ref{eq:IIy10}). We repeat the same procedure with
$\delta=0$ to arrive at the expression for $C(0)$. \hfill{} $\blacksquare$

\noindent \textbf{Proof of Lemma \ref{lemma:large-b}:} In the limit
$b\rightarrow\infty$, we can show that $\Delta_{IE}=\xi_{I}-\xi_{E}\rightarrow0$
and $\Delta_{E0}=\xi_{E}-\xi_{0}\rightarrow-\infty$ are the only
correct asymptotic behaviors. We notice that a necessary condition
for the firm at time $t$ to have non-negative return from investment
is that the boosted profit rate $X_{t}+b$ exceeds $\xi_{1}$, so
$\xi_{I}+b>\xi_{1}$ must be satisfied. Hence, in the limit $b\rightarrow\infty$,
$e^{\lambda(\Delta_{IE}+\Delta_{E0}+b+\xi_{0}-\xi_{1})}$ is bounded
by $1$ because $\Delta_{IE}+\Delta_{E0}+b+\xi_{0}-\xi_{1}=\xi_{I}+b-\xi_{1}>0$
and $\lambda<0$. 

From RHS of Eq. (\ref{eq:y11}), the leading-order term of $\Delta_{E0}$
is $-g$. We claim that the second-leading-order term of $\Delta_{E0}$
is a positive constant, independent of $g$. Suppose that the second-order
term of $\Delta_{E0}$ grows in $g$, but does so more slowly than
$g$. Then the first and second leading-order terms of Eqs. (\ref{eq:y11})
and (\ref{eq:y12}) are $-ge^{-\gamma_{p}\Delta_{IE}}=-g+g\gamma_{p}\Delta_{IE}(1+o(1))$
and $-ge^{-\gamma_{n}\Delta_{IE}}=-g+g\gamma_{n}\Delta_{IE}(1+o(1))$
respectively, which are inconsistent because $\gamma_{p}\neq\gamma_{n}$.
Thus, the second leading order term of $\Delta_{E0}$ is a constant
independent of $g$. Hence, we can express $\Delta_{E0}$ as in Eq.
(\ref{eq:y1-large-eta}) where $\theta$ is a constant yet to be determined.
Then Eqs. (\ref{eq:y11}) and (\ref{eq:y12}) can be re-expressed
as 
\begin{eqnarray}
\Delta_{E0} & = & -g+g\gamma_{p}\Delta_{IE}+(\lambda^{-1}-\gamma_{n}^{-1})e^{\lambda(\theta-\delta/\alpha+k\alpha+\xi_{0}-\xi_{1})}+o(1)\:,\label{eq:y1-eta-large-interim}\\
\Delta_{E0} & = & -g+g\gamma_{n}\Delta_{IE}+(\gamma_{p}^{-1}-\gamma_{n}^{-1})+(\lambda^{-1}-\gamma_{p}^{-1})e^{\lambda(\theta-\delta/\alpha+k\alpha+\xi_{0}-\xi_{1})}+o(1)\;.\label{eq:y12-eta-large-interim}
\end{eqnarray}
 Thus, the leading-order term of $\Delta_{IE}$ converges to zero
at least as fast as $g^{-1}$ because otherwise $\Delta_{E0}$ has
a second leading order term growing in $g$. Let us set $\Delta_{IE}=C/g+o(g^{-1})$
for some constant $C$. From Eqs. (\ref{eq:y1-large-eta}), (\ref{eq:y1-eta-large-interim})
and (\ref{eq:y12-eta-large-interim}), we arrive at $C=-(\gamma_{p}\gamma_{n})^{-1}(1-\lambda\theta-\lambda/\gamma_{n})$
where $\theta$ satisfies Eq. (\ref{eq:theta-eq.}). \hfill{}$\blacksquare$
\end{document}